\documentclass[12pt]{article}

\usepackage{array}
\usepackage{graphicx}
\usepackage{amsmath}
\usepackage{amsthm}
\usepackage{amssymb}
\usepackage{hhline}

\newcommand{\D}{\Delta}
\newcommand{\be}{\begin{equation}}
\newcommand{\ee}{\end{equation}}
\newcommand{\bsube}{\begin{subequations}}
\newcommand{\esube}{\end{subequations}}
\newcommand{\ba}{\begin{array}}
\newcommand{\ea}{\end{array}}
\newcommand{\To}{\rightarrow}
\newcommand{\bea}{\begin{eqnarray}}
\newcommand{\eea}{\end{eqnarray}}
\newcommand{\bc}{\begin{center}}
\newcommand{\ec}{\end{center}}

\newcommand{\dst}{\displaystyle}
\newcommand{\const}{{\rm const}}
\newcommand{\tu}{\tilde{u}}
\newcommand{\tv}{\tilde{v}}
\newcommand{\tw}{\tilde{w}}
\newcommand{\bu}{\bar{u}}
\newcommand{\dt}{\Delta t}
\newcommand{\dx}{\Delta x}
\newcommand{\F}{{\mathcal F}}
\newcommand{\sech}{{\rm sech}\,}
\newcommand{\usol}{u_{\rm sol}}
\newcommand{\Usol}{U_{\rm sol}}
\newcommand{\Ub}{U_{\rm b}}

\newcommand{\omsol}{\omega_{\rm sol}}
\newcommand{\ompw}{\omega_{\rm pw}}
\newcommand{\upw}{u_{\rm pw}}
\newcommand{\ub}{u_{\rm b}}

\title{\bf Instability of the finite-difference split-step
method on the background of
a soliton of the nonlinear Schr\"odinger equation}

\author{ T.I. Lakoba\thanks{{\tt tlakoba@uvm.edu}  } \\ \\
Department of Mathematics and Statistics,  University of Vermont, \\ Burlington, VT 05401, USA }

\begin{document}

\maketitle

\begin{abstract}
We consider the implementation of the split-step method where the 
linear part of the nonlinear Schr\"odinger equation is solved using a
finite-difference discretization of the spatial derivative.
The von Neumann analysis predicts that this method is unconditionally stable on 
the background of a
constant-amplitude plane wave. However, simulations show that the method can become unstable
on the background of a soliton. We present an analysis explaining this instability.
Both this analysis and the instability itself are substantially different from those
of the Fourier split-step method, which computes the spatial derivative by spectral discretization.
We also found that the modes responsible for the numerical instability are supported
by the sides of the soliton, in contrast to unstable modes of linearized nonlinear wave
equations, which (the modes) are supported by the soliton's core.
\end{abstract}

\bigskip

{\bf Keywords}: \ 
Finite-difference methods, Numerical instability, Nonlinear evolution equations.

\bigskip

{\bf PACS number(s)}: \ 02.60.-x, 02.60.Cb, 02.60.Lj, 02.70.Bf, 42.65.Tg.

\newpage

\section{Introduction}

The split-step method (SSM), also known as the operator-splitting method, is widely used
in numerical simulations of evolutionary equations that arise in diverse areas of science:
nonlinear waves, including nonlinear optics and Bose--Einstein condensation
\cite{HardinTappert_73}--\cite{Caliari2009}, atomic physics 
\cite{Bandrauk93,Bandrauk07}, studies of advection-reaction-diffusion equations
\cite{SSMReactDiff_99}--\cite{porousmed_2009}, and relativistic quantum mechanics
\cite{KGE_2012}. In this paper we focus on the SSM applied to the nonlinear
Schr\"odinger equation (NLS)
\be
i u_t - \beta u_{xx} + \gamma u |u|^2 =0, \qquad u(x,0)=u_0(x).
\label{e_01}
\ee
Although the real-valued constants $\beta$ and $\gamma$ in \eqref{e_01} can be scaled
out of the equation, we will keep them in order to distinguish the contributions of
the dispersive ($u_{xx}$) and nonlinear ($u|u|^2$) terms. Without loss of generality
we will consider $\gamma > 0$; then solitons of \eqref{e_01} exist for $\beta < 0$. 

The idea of the SSM is that \eqref{e_01} can be easily solved analytically when either
the dispersive or the nonlinear term is set to zero. This alows one to seek an approximate
numerical solution of \eqref{e_01} as a sequence of steps which alternatively account
for dispersion and nonlinearity:
\be
\ba{llr}
  \mbox{for $n$ from $1$ to $n_{\max}$ do:} & & \\
    & \hspace*{-4cm} \bar{u}(x) = u_n(x)\,\exp\big(i\gamma |u_n(x)|^2 \dt \big) 
    & \mbox{(nonlinear step)} \vspace{0.2cm}\\
    & \hspace*{-4cm} u_{n+1}(x) = 
      \left\{ \ba{l} \mbox{solution of \ $iu_t=\beta u_{xx}$ at $t=\dt$} \\
              \mbox{with initial condition $u(x,0)=\bar{u}(x)$}
      \ea  \right.
    & \mbox{(dispersive step)}
\ea
\label{e_02}
\ee
where the implementation of the dispersive step will be discussed below. 
In \eqref{e_02}, $\dt$ is the time step of the numerical integration and
$u_n(x) \equiv u(x,n\dt)$. Scheme \eqref{e_02} can yield a numerical solution
of \eqref{e_01} whose accuracy is $O(\dt)$. Higher-order schemes, yielding
more accurate solutions (e.g., with accuracy $O(\dt\,^2)$, $O(\dt\,^4)$, etc.),
are known \cite{Strang,Yevick91,Bandrauk93}, but here we will restrict our
attention to the lowest-order scheme \eqref{e_02}; see also the 
paragraph after Eq.~\eqref{e_33} below.

The implementation of the dispersive step in \eqref{e_02} depends on the numerical
method by which the spatial derivative is computed. In most applications, it is
computed by the Fourier spectral method:
\be
u_{n+1}(x) = \F^{-1} \left[ \,\exp(i\beta k^2 \dt) \; \F[ \bu(x)]\;\right]\,.
\label{e_03}
\ee
Here $\F$ and $\F^{-1}$ are the discrete Fourier transform and its inverse,
$k$ is the discrete wavenumber:
\be
-\pi/\dx \le k \le \pi/\dx,
\label{e_04}
\ee
and $\dx$ is the mesh size in $x$. However, the spatial derivative in \eqref{e_02}
can also be computed by a finite-difference (as opposed to spectral) method
\cite{WH,SelAreasCommun_1997,JCP1999,Faou_2011}.
For example, 
using the central-difference discretization of $u_{xx}$ and the Crank--Nicolson
method, the dispersion step yields:
\be
i \frac{u_{n+1}^m - \bu^m}{\dt} \,=\, \frac{\beta}2 \,
\left( \frac{u_{n+1}^{m+1}-2u_{n+1}^m+u_{n+1}^{m-1}}{\dx\,^2} + 
\frac{\bu^{m+1}-2\bu^m+\bu^{m-1}}{\dx\,^2} \right),
\label{e_05}
\ee
where $u_n^m \equiv u(x_m,n\dt)$, $x_m$ is a point in the discretized spatial
domain: \ $-L/2 < x_m < L/2$, and $L$ is the the length of the domain. 
Recently, solving the dispersive step of \eqref{e_02} by a finite-difference 
method has found an application in the electronic post-processing of the optical
signal in fiber telecommunications \cite{OE_2008}. That post-processing must be
done in real time, and hence the speed of its implementation becomes a key factor.
Ref.~\cite{JLT_2005} showed that finite-difference implementations (referred to there
as the Finite and Infinite Impulse Response filters; see its Eqs. (13) and (22))
of the dispersive step of \eqref{e_02} provide considerable speed improvement
over the Fourier-based solution \eqref{e_03};
see also Sec.~2.4.2 in \cite{Agrawal_book}.
Finally, let us note that the version of the NLS where the second derivative is 
replaced by its
finite-difference approximation, as on the right-hand side (r.h.s.) of \eqref{e_05},
is of considerable interest also in its own right \cite{DNLS}.

In what follows we assume periodic boundary conditions:
\be
u(-L/2,t)=u(L/2,t);
\label{e_06}
\ee
the case of other types of boundary conditions is considered in Appendix A.
We will refer to the SSM with spectral \eqref{e_03} and finite-difference 
\eqref{e_05} implementations of the dispersive step in \eqref{e_02} as
s-SSM and fd-SSM, respectively. 
Our focus in this paper will on the fd-SSM.

\begin{figure}[h]
\vspace{-1.6cm}
\mbox{ 
\begin{minipage}{7cm}
\rotatebox{0}{\resizebox{7cm}{9cm}{\includegraphics[0in,0.5in]
 [8in,10.5in]{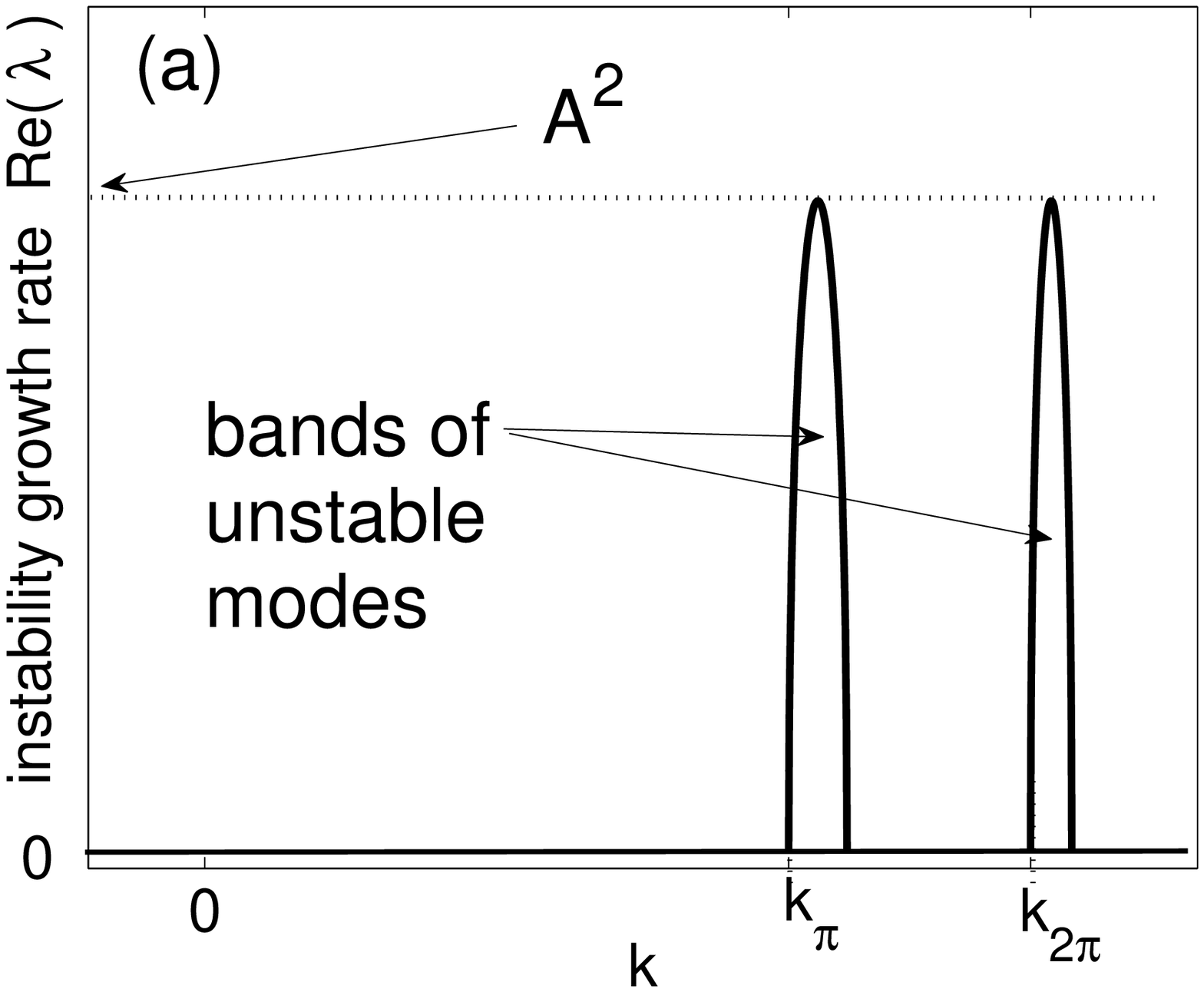}}}
\end{minipage}
\hspace{0.1cm}
\begin{minipage}{7cm}
\rotatebox{0}{\resizebox{7cm}{9cm}{\includegraphics[0in,0.5in]
 [8in,10.5in]{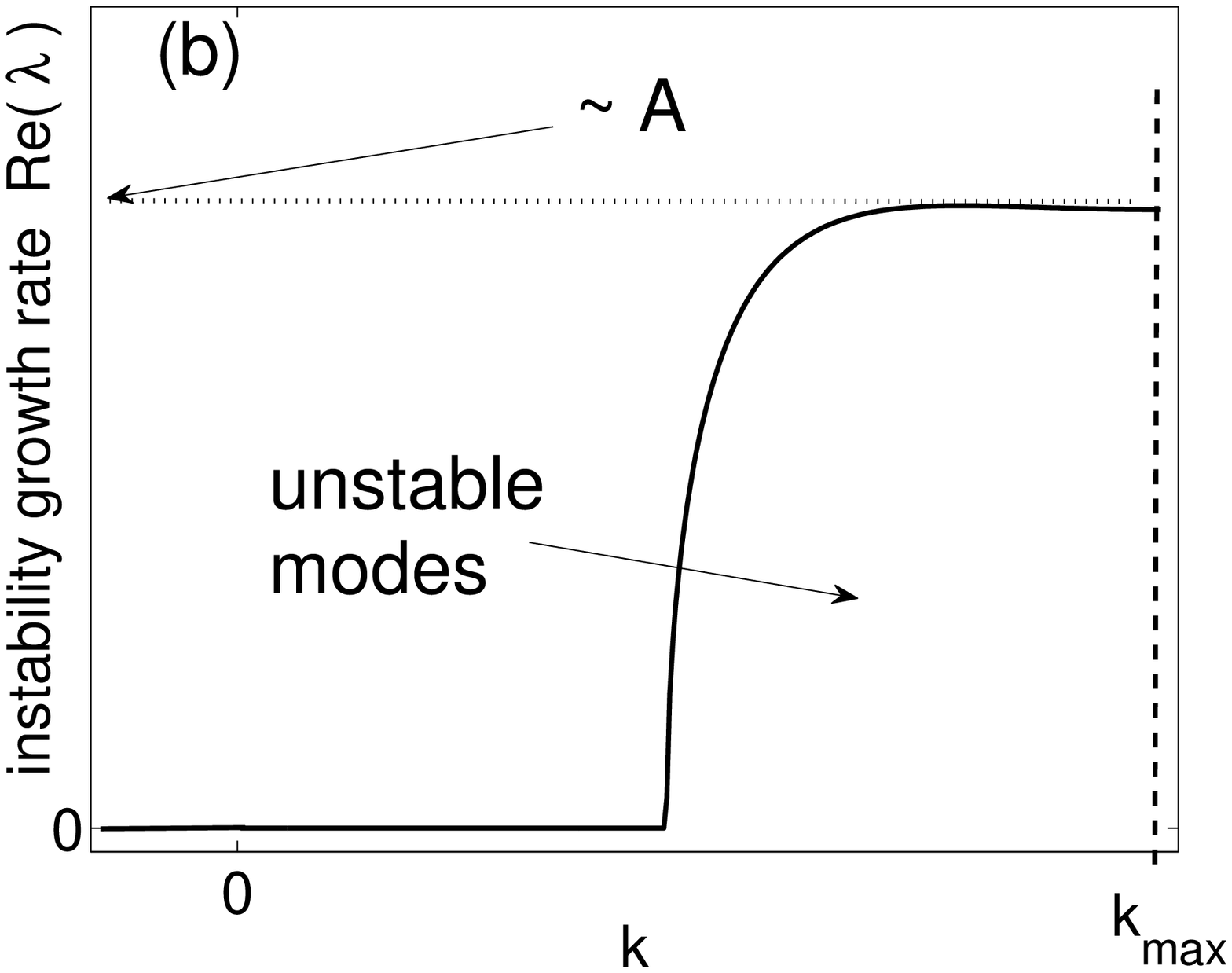}}}
\end{minipage}
 }
\vspace{-1.6cm}
\caption{Growth rate of numerical instability of the s-SSM (a) and
fd-SSM (b) on the plane-wave background. 
The dotted horizontal line indicates how the maximum growth rate depends on the
wave's amplitude. In (a), $k_{m\pi}$, $m=1,2,\ldots$ are the wavenumbers where the
$m$th resonance condition holds (see \cite{ja}): $|\beta|k_{m\pi}^2\dt = m\pi$. 
}
\label{fig_1}
\end{figure}

Weideman and Herbst \cite{WH} used the von Neumann analysis to show that 
both versions, s- and fd-, of the SSM can become unstable when the background
solution of the NLS is a plane wave:
\be
\upw = (A/\sqrt{\gamma}) \; e^{i\omega_{\rm pw}t}, \qquad A=\const, \quad
\omega_{\rm pw} = |A|^2.
\label{e_07}
\ee
Specifically, they linearized the SSM equations on the background of \eqref{e_07}:
\be
u_n=\upw + \tu_n, \qquad |\tu_n| \ll |u_n|
\label{e_08}
\ee
and sought the numerical error in the form
\be 
\tu_n = \tilde{A}\,e^{\lambda t_n - ikx}, \qquad \tilde{A}=\const.
\label{e_09}
\ee
The SSM is said to be unstable when for a certain wavenumber $k$ one has: \ 
(i) \ Re$(\lambda) > 0$ in \eqref{e_09}, but \ (ii) \ the corresponding Fourier
mode in the original equation \eqref{e_01} is linearly stable. Weideman and Herbst
found that the s- and fd-SSMs on the background \eqref{e_07} become unstable when
the step size $\dt$ exceeds:
\be
\dt_{\rm thr,\,s} \approx \dx\,^2/(\pi |\beta|), \qquad 
\mbox{for s-SSM \eqref{e_02} \& \eqref{e_03}}
\label{e_10}
\ee
and
\be
\dt_{\rm thr,fd}=\dx/\sqrt{2|\beta|\,|A|^2} \quad \underline{\mbox{only for $\beta>0$}}, \qquad 
\mbox{for fd-SSM \eqref{e_02} \& \eqref{e_05}}
\label{e_11}
\ee
respectively. 
Note that for $\beta<0$, the fd-SSM simulating a solution close to the plane wave 
\eqref{e_07} is unconditionally stable. 
Typical dependences of the instability growth rate, Re$(\lambda)>0$,
on the wavenumber is shown in Fig.~\ref{fig_1}. 
Let us emphasize that the SSM is unstable for $\dt>\dt_{\rm thr}$ even
though both its constituent steps, \eqref{e_02} and either \eqref{e_03} or
\eqref{e_05}, are numerically stable for any $\dt$. 

Solutions of the NLS (and of other evolution equations) that are of
practical interest are considerably more complicated than a constant-amplitude
solution \eqref{e_07}. To analyze stability of a numerical method that is
being used to simulate a {\em spatially varying} solution, one often
employs the so-called ``principle of frozen coefficients" \cite{vonNeumannRichtmeyer}
(see also, e.g., \cite{Trefethen_book,KGE_2012}). According to that principle,
one assumes some constant value for the solution $u$ and then linearizes the
equations of the numerical method to determine the evolution of the numerical error
(see \eqref{e_08} and \eqref{e_09}). However, as we show below, this principle
applied to the SSM fails to predict, even qualitatively, important features
of the numerical instability.

As a first step towards developing a general framework of stability analysis on the
background of a spatially varying solution, we analyzed \cite{ja} the instability
of the s-SSM on the background of a soliton of the NLS:
\be
\usol(x,t) = A \sqrt{2/\gamma} \; \sech(Ax/\sqrt{-\beta})\; e^{i\omega_{\rm sol}t}
 \,\equiv \, \Usol(x)\,e^{i\omega_{\rm sol}t}, 
 \qquad \omega_{\rm sol}=A^2.
\label{e_12}
\ee
First, we demonstrated numerically that the instability growth rate in this case
is very sensitive to the time step $\dt$ and the length $L$ of the spatial domain;
also, its dependence on the wavenumber is quite different from that shown in
Fig.~\ref{fig_1}(a). Moreover, the instability on the background of, say, two 
well-separated (and hence non-interacting) solitons can be completely different from 
that on the background of one of these solitons. To our knowledge, such features
of the numerical instability had not been reported for other numerical methods. Moreover,
they could not be predicted based on the principle of frozen coefficients. We then
demonstrated that all those features could be explained by analyzing an equation
satisfied by the numerical error of the s-SSM:
\be
i\tv_t - \omsol \tv -\beta(\tv_{xx}+k^2_{\pi}\tv) + \gamma |\usol|^2 
(2\tv + \tv^*) = 0,
\label{e_13}
\ee
where $\tv(x,t)$ is proportional to the continuous counterpart of 
$\tu_n(x)\equiv \tu(x,n\dt)$ defined in \eqref{e_08}, and $k_{\pi}$ is defined
in the caption to Fig.~\ref{fig_1}. Note that \eqref{e_13} is similar, but not 
equivalent, to the NLS linearized about the soliton:
\be
i\tu_t - \omsol \tu -\beta\tu_{xx} + \gamma |\usol|^2 
(2\tu + \tu^*) = 0.
\label{e_14}
\ee
The extra $k^2_{\pi}$-term in \eqref{e_13} indicates that the potentially 
unstable numerical error of the s-SSM has a wavenumber close to $k_{\pi}$.

In this paper we present an analysis of the numerical instability of the
fd- (as opposed to s-) SSM on the soliton background. This analysis is based
on an equation for the numerical error which, as \eqref{e_13}, is a modified 
form of the linearized NLS. However, both that equation and its analysis are
{\em qualitatively} different from those \cite{ja} for the s-SSM. Moreover,
the modes that render the s- and fd-SSMs unstable are also qualitatively
different. Namely, for the s-SSM, the numerically unstable modes contain just a few
Fourier harmonics and hence are not spatially localized. On the contrary, 
the modes making the fd-SSM unstable are localized and are supported by the
sides (i.e., tails) of the soliton. To our knowledge, such ``tail-supported"
localized modes, as opposed to those supported by the soliton's core, 
have not been reported before.

The main part of this manuscript is organized as follows. In Sec.~II we present
simulation results showing the development of numerical instability of the fd-SSM
applied to a soliton. In Sec.~III we derive an equation (a counterpart of 
\eqref{e_13}) governing the evolution of the numerical error, and in Sec.~IV
obtain its localized solutions that grow exponentially in time. 
In Sec.~V we summarize the results. In Appendix A we show how our analysis
can be modified for boundary conditions other than periodic. In Appendix B
we describe the numerical method used to obtain the localized solutions in Sec.~IV.
In Appendix C we discuss how the instability sets in.


\section{Numerics of fd-SSM with soliton background}

We numerically simulated Eq. \eqref{e_01} with $\beta = -1$, $\gamma=2$, and the periodic
boundary conditions \eqref{e_06} via the fd-SSM algorith \eqref{e_02} \& \eqref{e_05}.
The initial condition was the soliton \eqref{e_12} with $A=1$:
\be 
u_0(x) = \sech(x) + \xi(x);
\label{e_15}
\ee
the noise component $\xi(x)$ with zero mean and the standard deviation $10^{-10}$
was added in order to reveal the unstable Fourier components sooner than if they 
had developed from a round-off error. 

Below we report results for two values of the spatial mesh size $\dx=L/N$, where
$N$ is the number of grid points: \ $\dx=40/2^9$ and $\dx=40/2^{10}$. We verified that,
for a fixed $\dx$, the results are insensitive to the domain's length $L$
(unlike they are for the s-SSM \cite{ja}) as long as $L$ is sufficiently large. 
Also, at least within the range of $\dx$ values considered, the results depend
on $\dx$ monotonically (again, unlike for the s-SSM).

\begin{figure}[h]
\vspace{-1.6cm}
\mbox{ 
\begin{minipage}{7cm}
\rotatebox{0}{\resizebox{7cm}{9cm}{\includegraphics[0in,0.5in]
 [8in,10.5in]{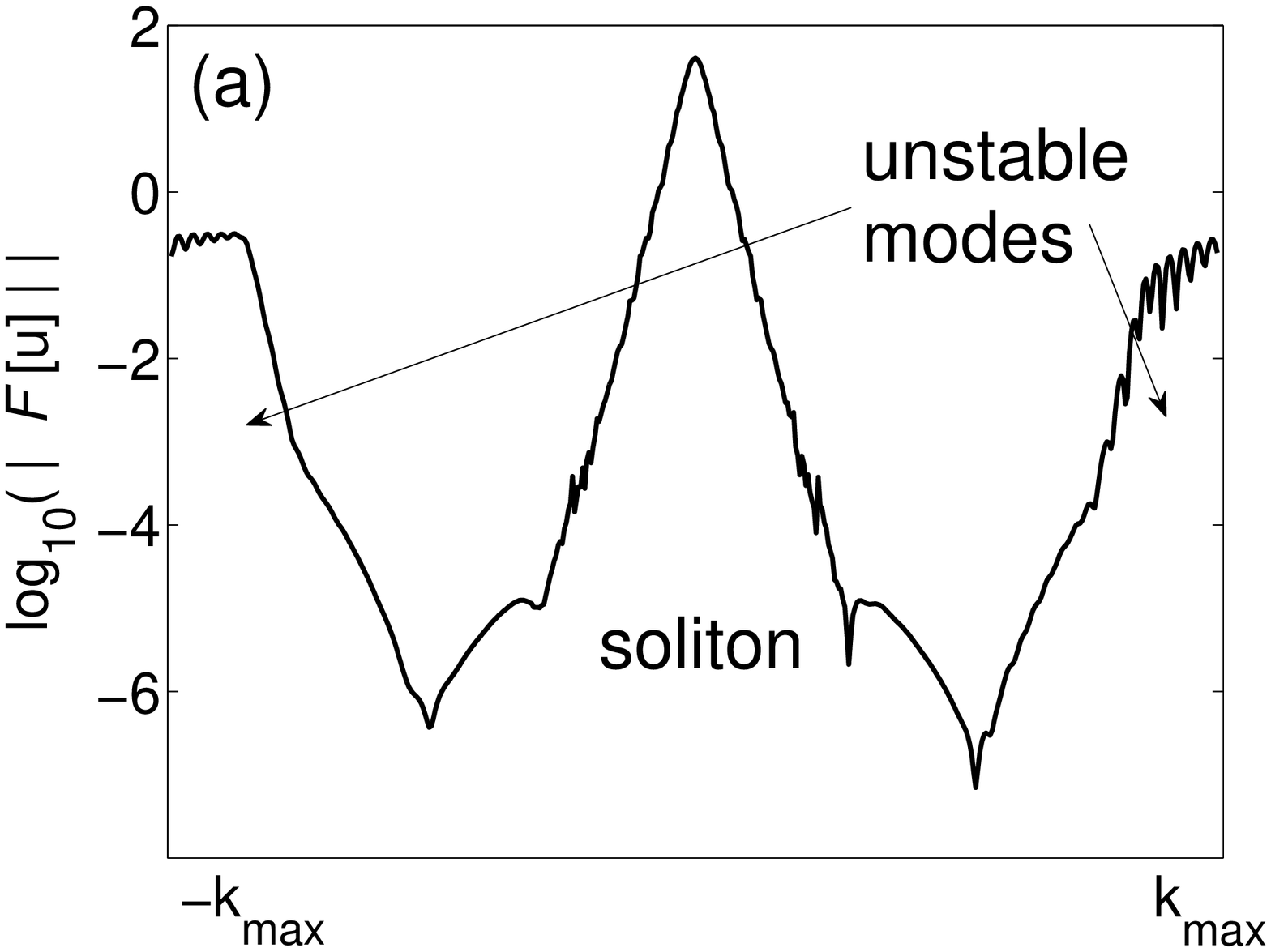}}}
\end{minipage}
\hspace{0.1cm}
\begin{minipage}{7cm}
\rotatebox{0}{\resizebox{7cm}{9cm}{\includegraphics[0in,0.5in]
 [8in,10.5in]{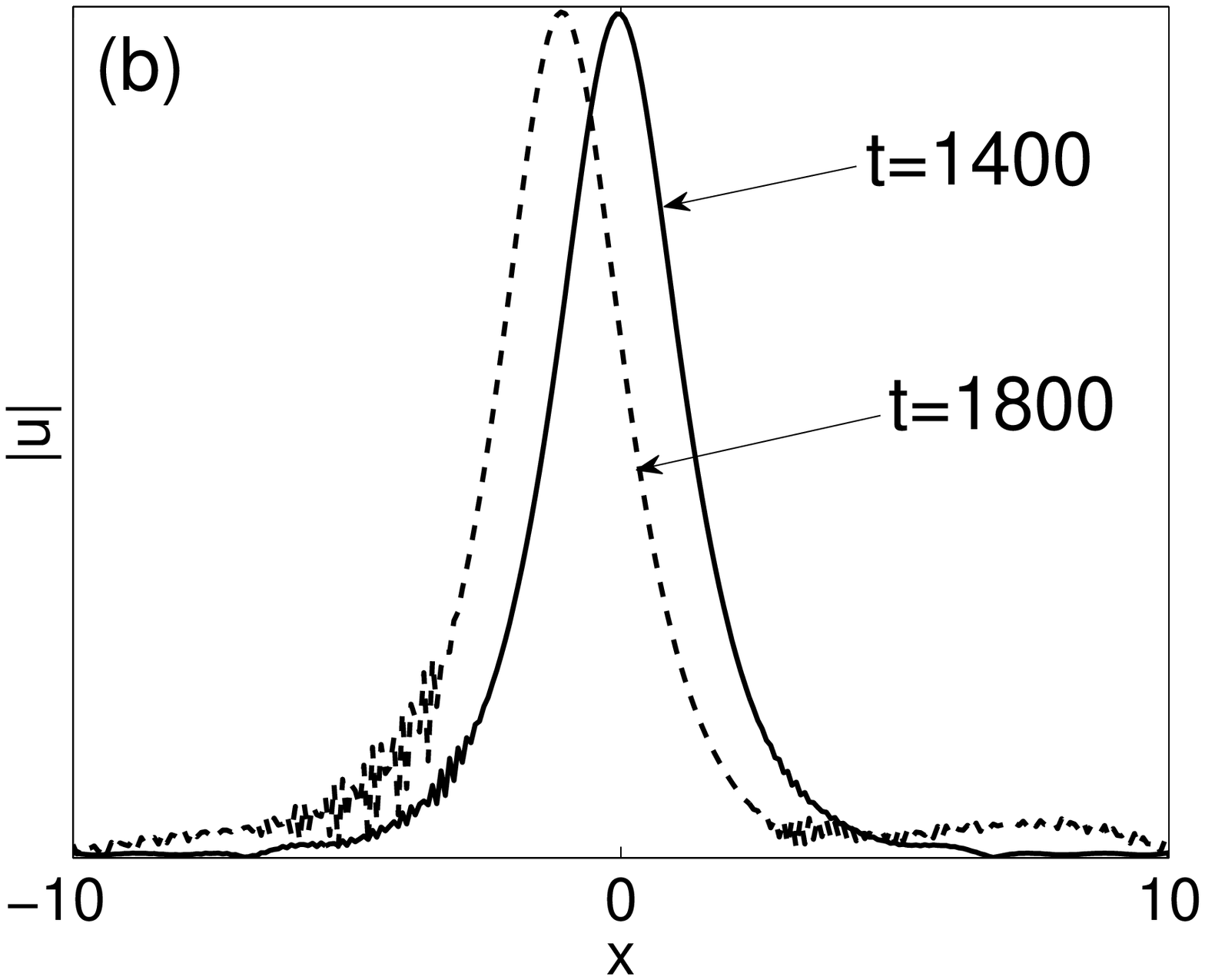}}}
\end{minipage}
 }
\vspace{-1.6cm}
\caption{ See explanation in the text.
}
\label{fig_2}
\end{figure}

First, let us remind the reader that the analysis of \cite{WH} on a constant-amplitude
background \eqref{e_07} for $\beta <0$ predicted that the fd-SSM should be stable for 
any $\dt$.\footnote{
Recall that the stability of instability of the numerical method is in no way related
to that of the actual solution. In fact, the plane wave \eqref{e_07} is well-known
to be modulationally unstable for $\beta<0$, while it is modulationally stable 
for $\beta > 0$.}
For the soliton initial condition \eqref{e_15} and the parameters stated above, 
our simulations showed that the 
numerical solution becomes unstable for $\dt> \dx$. For future use we introduce a notation:
\be
C = (\dt/\dx)^2\,.
\label{e_16}
\ee
In Fig.~\ref{fig_2}(a) we show the Fourier spectrum of the numerical solution of 
\eqref{e_01}, \eqref{e_15} obtained by the fd-SSM with $C=1.05$ 
(i.e., slightly above the instability threshold)  at $t=1400$.
The numerically unstable modes are seen around the end points of the spectral axis.
At $t=1400$, these modes are still small enough so as not to cause 
visible damage to the soliton: see the solid curve in Fig.~\ref{fig_2}(b).
However, at a later time, the soliton begins to drift: see the dashed line in
Fig.~\ref{fig_2}(b), that shows the numerical solution at $t=1800$. Such a drift 
may persist over a long time: e.g., for $C=1.05$, the soliton 
still keeps on moving at $t\sim 4000$. However, eventually it gets overcome by noise and 
loses its identity. 

We observed the same scenario for several different values of $\dx$, $L$, and $C$
(for $C>1$). The direction of the soliton's drift appears to be determined by the
initial noise; this direction is {\em not} affected by the placement of the initial
soliton inside the spatial domain. The time when the drift's onset becomes visible 
decreases, and the drift's velocity increases, as $C$ increases.

\begin{figure}[h]
\vspace{-1.6cm}
\centerline{ 
\begin{minipage}{7cm}
\rotatebox{0}{\resizebox{7cm}{9cm}{\includegraphics[0in,0.5in]
 [8in,10.5in]{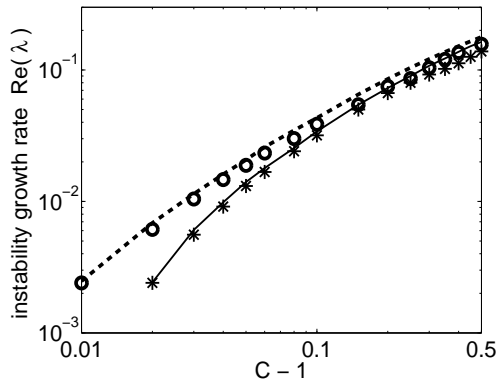}}}
\end{minipage}
 }
\vspace{-1.6cm}
\caption{Growth rate of the numerical instability for $\dx=40/2^9$ (solid line --- 
analysis of Sec.~IV, stars --- computed by \eqref{e_17}) and for 
$\dx=40/2^{10}$ (dashed line --- 
analysis of Sec.~IV, circles --- computed by \eqref{e_17}).
}
\label{fig_3}
\end{figure}

The soliton's drift is a nonlinear stage of the development of the numerical 
instability. It will be explained in Sec.~V. In the linear stage, the numerically
unstable modes are still small enough so that they do not visibly affect the soliton
or one another. To describe this stage, we computed a numeric approximation to the
instability growth rate \ Re$(\lambda)$ defined in \eqref{e_09}:
\be
{\rm Re}\,(\lambda)|_{\rm computed} = 
\frac{ \ln\left( \ba{c} \max \left| \F[u](k)\right| \;\; \mbox{for} \\ 
                        \mbox{$k \sim k_{\max}$ at time$=t$} \ea \right) -
       \ln \left( \ba{c} \mbox{noise floor} \\ \mbox{at time$=0$} \ea \right) }{t},
\label{e_17}
\ee
where $k_{\max}=\pi/\dx$ (see \eqref{e_04}).
The so computed
values of the instability growth rate are shown in Fig.~\ref{fig_3} along with the results
of a semi-analytical calculation presented in Sec.~IV. 

The above numerical results motivate the following three questions: \ (i) explain the observed
instability threshold $\dt$ (see the sentence before \eqref{e_16}); \ (ii) identify the modes
responsible for the numerical instability; and \ (iii) calculate the instability growth rate 
(see Fig.~\ref{fig_3}). In Sec.~IV we will give an approximate analytical answer to question (i).
However, answers to questions (ii) and (iii) will be obtained only semi-analytically, i.e., via
numerical solution of a certain eigenvalue problem.


\section{Derivation of equation for numerical error for fd-SSM}

\begin{figure}[h]
\vspace{-1.6cm}
\centerline{ 
\begin{minipage}{7cm}
\rotatebox{0}{\resizebox{7cm}{9cm}{\includegraphics[0in,0.5in]
 [8in,10.5in]{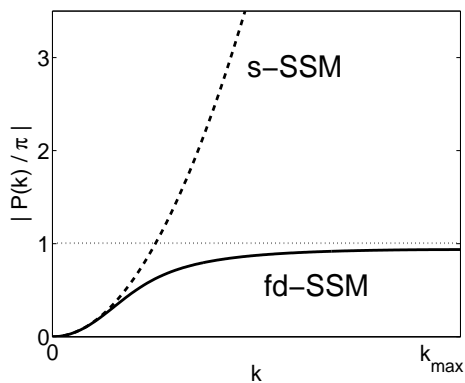}}}
\end{minipage}
 }
\vspace{-1.6cm}
\caption{ Normalized phase: \ $|\beta|k^2\dt$ for the s-SSM (dashed) and as given by \eqref{e_19} 
for the fd-SSM (solid). In both cases, $r=5$. The horizontal line indicates the condition
of the first resonance: $|P(k)|=\pi$.
}
\label{fig_4}
\end{figure}

Here we will derive a modified linearized NLS that is a counterpart of \eqref{e_13},
which was obtained for the s-SSM in \cite{ja}. In view of periodic boundary conditions
\eqref{e_06}, the finite-difference implementation \eqref{e_05} of the dispersive step 
in \eqref{e_02} can be written as
\be
u_{n+1}(x) = \F^{-1} \left[ e^{iP(k)}\,\F\left[ \bu(x) \right] \; \right],
\label{e_18}
\ee
\be
e^{iP(k)} \equiv \frac{1+2i\beta r \sin^2(k\dx/2) }{1-2i\beta r \sin^2(k\dx/2) }
 = \exp\left[ 2i{\rm arctan}\big( 2\beta r \sin^2(k\dx/2)\,\big) \right], 
 \qquad r = \frac{\dt}{\dx\,^2},
\label{e_19}
\ee
where $\F$, $\F^{-1}$ were defined after \eqref{e_03}. For $|k\dx|\ll 1$, the 
exponent in \eqref{e_19} equals that in \eqref{e_03}. However, for $|k\dx|>1$, they
differ substantially: see Fig.~\ref{fig_4}.
It is this difference that leads to the instabilities of the s- and fd-SSMs being
qualitatively different.

Using Eqs. \eqref{e_02} and \eqref{e_18}, one can write, similarly to Eq. (3.1) in
\cite{ja}, a linear equation satisfied by a small numerical error $\tu_n$ of the 
fd-SSM:
\be
\F \left[ \tu_{n+1} \right] = e^{iP(k)} \F
 \left[ e^{i\gamma |\ub|^2\dt} \big( \tu_n + 
        i\gamma \dt ( \ub^2 \tu_n^* + |\ub|^2 \tu_n )\, \big) \; \right]\,.
 \label{e_20}
 \ee
 Here $\tu_n$ is defined similarly to \eqref{e_08}, with $\ub$ being either
 $\upw$ or $\usol$, depending on the background solution. The exponential growth
 of $\tu_n$ can occur only if there is sufficiently strong coupling between
 $\tu_n$ and $\tu_n^*$ in \eqref{e_20}. This coupling is the strongest when
 the temporal rate of change of the realtive phase between these two terms is
 the smallest. In \cite{ja} we showed that this rate can be small only for those $k$
 where the exponent $P(k)$ is close to a multiple of $\pi$. Using \eqref{e_19}
 (see also Fig.~\ref{fig_4}),
 we see that this can occur only for sufficiently high $k$ where 
 $\sin^2(k\dx/2)=O(1)$ rather than $O(\dx\,^2)$. Then:
 \bea
 P(k) & = & \pi - \frac1{\beta r \sin^2(k\dx/2)} + O\left( \frac1{r^3}\right)
  \nonumber  \\
  & = & \pi - \frac1{\beta r} - \frac{(k-k_{\max})^2 \dx\,^2}{4\beta r} + 
        O\left( \frac1{r^3} + \frac{\big( (k-k_{\max}) \dx\,\big)^4}{r} \right),
 \label{e_21}
 \eea
 where $k_{\max}=\pi/\dx$. We have also used that 
 \be
 r = \dt/\dx\,^2 = C/\dt \gg 1,
 \label{e_22}
 \ee
 given that the numerical instability was observed in Sec.~II for $C=O(1)$. 
 
 In order for the truncation of the Taylor expansions used in \eqref{e_21} to be
 self-consistent, one needs to assume that
 $$
 O(1/r^2) < O\big( (k-k_{\max})^2\dx\,^2 \big) < O(1);
 $$
 note that $O(1/r^2) = O(\dx\,^2)$. We arbitrarily declare $(k-k_{\max})^2\dx\,^2$
 to be in the middle of that range: $(k-k_{\max})^2\dx\,^2 = O(\dx)$, which implies that
 the range of wavenumbers for which the instability develops satisfies:
 \be
 |k-k_{\max}| = O(1/\sqrt{\dx}) \ll k_{\max} = \pi/\dx.
 \label{e_23}
 \ee
 This turns out to be consistent with the result shown in Fig.~\ref{fig_2}(a). 
 
 Substituting the first three terms on the r.h.s. of \eqref{e_21} into \eqref{e_20}, using 
 \eqref{e_22}, and introducing a new variable 
 \be
 \tv_n = \left( e^{i\pi} \right)^n \tu_n = (-1)^n \tu_n,
 \label{e_24}
 \ee
 one obtains:
 \bea
 \F [\tv_{n+1}] & = & \exp\left( -\frac{i\dt}{C \beta} 
                       \left\{1 + \frac{ (k-k_{\max})^2\dx\,^2}4 \right\} \right) \,\times
  \nonumber \\
   &  & \F \left[ e^{i\gamma |\ub|^2\dt}\, 
          \left\{ \tv_n + i\gamma \dt ( \ub^2\tv_n^* + |\ub|^2 \tv_n ) \right\} \,
           \right].
 \label{e_25}
 \eea
 Note that \eqref{e_25} describes a {\em small} change of $\tv_n$ occurring over the step
 $\dt$, because for $\dt\To 0$, the r.h.s. of that equation reduces to 
 $\F[\tv_n]$. Therefore we can approximate the {\em difference} equation \eqref{e_25}
 by a {\em differential} equation. To that end, first recall from \eqref{e_23}
 that the wavenumbers of $\tv_n$ are on the order of $k_{\max}$; hence we seek
 \be
 \tv_n(x) = e^{i k_{\max} x}\, \tw_n(x).
 \label{e_26}
 \ee
 The effective wavenumber of $\tw_n$ is then $(k-k_{\max})$, and according to
 \eqref{e_23} $\tw_n$ varies slowly over the scale $O(\dx)$. Introducing the scaled
 space and wavenumber by
 \be
 \chi= x/\epsilon, \quad \kappa = (k-k_{\max})\epsilon, \qquad \epsilon = \dx/2,
 \label{e_27}
 \ee
 one rewrites \eqref{e_25} as:
 \be
 \F_{\rm sc} [\tw_{n+1}] = \exp\left( -\frac{i\dt}{C \beta} 
                        \{1 + \kappa^2 \} \right) \;
 \F_{\rm sc} \left[ e^{i\gamma |\ub|^2\dt}\, 
           \left\{ \tw_n + i\gamma \dt ( \ub^2\tw_n^* + |\ub|^2 \tw_n ) \right\} \,
            \right],
  \label{e_28}
  \ee
 where now $\F_{\rm sc}$ is the Fourier transform with respect to the scaled variables
 \eqref{e_27}. In handling the $\tv_n^*$ term in \eqref{e_25}, we have used the fact that
 on the spatial grid $x_m=m\dx$, one has:
 $$
 \tv_n^*(x_m) = e^{-i k_{\max}x_m} \tw_n^*(x_m) = e^{-i\pi m} \tw_n^*(x_m) 
 = e^{i\pi m} \tw_n^*(x_m) = e^{ik_{\max}x_m} \tw_n^*(x_m).
 $$

Next, the s-SSM \eqref{e_02}, \eqref{e_03} can be written as 
\be
\F [u_{n+1}] = e^{i\beta k^2 \dt}\; \F
 \left[ e^{i\gamma |u|^2 \dt}\, u \right].
\label{e_29}
\ee
When $|\beta|k^2\dt \ll 1$ and $\gamma |u|^2\dt \ll 1$, 
this is equivalent to the NLS \eqref{e_01}
{\em plus} a term proportional to 
\be
\dt\; \left[ \beta \partial^2_x, \; \gamma |u|^2 \right]_{-} u + O(\dt\,^2),
\label{e_30}
\ee
where $[\ldots,\; \ldots]_{-}$ denotes a commutator (see, e.g., Sec. 2.4 in 
\cite{Agrawal_book}). Equation \eqref{e_28} has the form of a linearized Eq. 
\eqref{e_29} with a different coefficient in the dispersion term and 
with an extra
phase. Therefore, \eqref{e_28} must be equivalent to a modified linearized NLS,
with the modification affecting only the corresponding terms:
\be
i\tw_t + (\tw_{\chi\chi}- \tw)/(C\beta) + \gamma (\ub^2 \tw^* + 2|\ub|^2 \tw) =0,
\label{e_31}
\ee
{\em plus} a term proportional to the linearized form of the commutator \eqref{e_30}.
Neglecting that latter term as small (of order $O(\dt)$) compared to the rest of
the expression and denoting \ $\psi=\tw\;\exp(-i\omega_{\rm b}t)$, we rewrite 
\eqref{e_31} as:
\be
i\psi_t + \delta\psi +\psi_{\chi\chi}/(C\beta) + \gamma \Ub^2(\epsilon\chi)\;
 (2\psi + \psi^*) =0,
\label{e_32}
\ee
where
\be
\delta = -\omega_{\rm b} - 1/(C\beta).
\label{e_33}
\ee
Here $\omega_{\rm b}$ is either $\ompw$ or $\omsol$, and $\Ub$ is either constant
or $\Usol$, depending on whether the background solution is a plane wave \eqref{e_07}
or a soliton \eqref{e_12}. The modified linearized NLS \eqref{e_32} for the fd-SSM
is the counterpart of Eq.~\eqref{e_13} that was derived for the s-SSM.

Our subsequent analysis of the instability of the first-order accurate
fd-SSM \eqref{e_02} \& \eqref{e_05}
will be based on Eq.~\eqref{e_32}. 
The instability of the second-order accurate version of this method,
where the order of the nonlinear and dispersive steps is alternated
in any two consecutive full time steps \cite{Strang}, is the same as 
that of the first-order version. The instability of higher-order 
versions (e.g., $O(\dt\,^4)$-accurate) can be studied similarly to how
that was done in Ref.~\cite{ja} for the s-SSM.

The boundary conditions satisfied by $\psi$ are still periodic:
\be
\psi(-L/(2\epsilon),\,t) = \psi(L/(2\epsilon),\,t).
\label{e_34}
\ee
This follows from the fact that $\tu_n(x)$ satisfies the periodic boundary conditions
\eqref{e_06} and from \eqref{e_26}, given that for $k_{\max}=\pi/\dx$ and 
$L/2=M\dx$ with some integer $M$,
$$
e^{-ik_{\max}L/2} = e^{-iM\pi} = e^{iM\pi} = e^{ik_{\max}L/2} .
$$

There are three differences between Eq.~\eqref{e_32} and the linearized NLS \eqref{e_14}.
Most importantly, \eqref{e_32} has the opposite sign of the dispersion term. This is
explained by the shape of the curve $P(k)$ 
for the fd-SSM in Fig.~\ref{fig_4} at high wavenumbers,
where the curvature is opposite to that at $k\approx 0$. 
Secondly, unlike the $(-\omsol)$-term in \eqref{e_14}, the $\delta$-term in \eqref{e_32}
with $\beta<0$ can be either positive or negative, depending on the value of $C$.
Thirdly, the ``potential" $\Ub^2(\epsilon\chi)$ (when $\Ub\equiv \Usol$) is a {\em slow}
function of the scaled variable $\chi$. That is, solutions of \eqref{e_32} that vary on 
the scale $\chi=O(1)$ ``see" the soliton as being very wide. This should also be
contrasted with the situation for the s-SSM, where the modes described by Eq.~\eqref{e_13}
``see" the soliton as being very narrow \cite{ja}.

Before proceeding to find unstable modes of Eq.~\eqref{e_32} with $\Ub\equiv \Usol$, let
us note that \eqref{e_32} with $\Ub=\const$ confirms the result of Ref.~\cite{WH}
regarding the instability of the fd-SSM on the plane-wave background. Namely, for
$\beta<0$, Eq.~\eqref{e_32} with $\Ub=\const$ describes the evolution of a small 
perturbation to the plane wave in the {\em modulationally  stable} case (see, e.g.,
Sec.~5.1 in \cite{Agrawal_book}). That is, for $\beta<0$, there is no numerical instability,
in agreement with \cite{WH}. On the other hand, for $\beta>0$, Eq.~\eqref{e_32}
describes the evolution of a small perturbation in the {\em modulationally  unstable}
case, and hence the plane wave of the NLS \eqref{e_01} can become numerically unstable.
The corresponding instability growth rate found from \eqref{e_32} and Eq. (5.1.8)
of \cite{Agrawal_book} can be shown to agree with the one that can be obtained from
Eq.~(37) and the next two unnumbered relations in \cite{WH}. An example of this growth
rate is shown in Fig.~\ref{fig_1}(b). Also, using our \eqref{e_32} 
and Eq.~(5.1.8) of Ref.~\cite{Agrawal_book}, 
the threshold value of $\dt$ can be shown to be given by \eqref{e_11},
in agreement with \cite{WH}.


\section{Unstable modes of modified linearized NLS for fd-SSM}

In this section we focus on the case where the background solution is a soliton;
hence $\beta < 0$ and $\Ub\equiv \Usol$ (see \eqref{e_12}). Substituting into
\eqref{e_32} and its complex conjugate the standard ansatz \cite{Kaup90}
\ $(\psi(\chi,t),\,\psi^*(\chi,t))=(\phi_1(\chi),\,\phi_2(\chi))\,e^{\lambda t}$ \ 
and using yet another rescaling:
\be
\ba{c}
\dst
X=\frac{A}{\sqrt{-\beta}}\chi \equiv \frac{2A}{\sqrt{-\beta}}\frac{x}{\dx},
\qquad
D = -\frac{C\beta^2}{A^2}\delta \equiv \beta^2\left(\frac1{\beta A^2} + C\right),
\vspace{0.2cm} \\
\dst
\Lambda = \frac{C\beta^2}{A^2}\lambda,  \qquad
V(y) = 2C\beta^2 \sech^2(y),
\ea
\label{e_35}
\ee
one obtains:
\be
\left( \partial_X^2 + D - V(\epsilon X) 
 \left( \ba{cc} 2 & 1 \\ 1 & 2 \ea \right) \,\right)
\vec{\phi} \,=\, i\Lambda \sigma_3 \vec{\phi},
\label{e_36}
\ee
where $\sigma_3={\rm diag}(1,-1)$ is a Pauli matrix, $\vec{\phi}=(\phi_1,\,\phi_2)^T$,
and $T$ stands for a transpose. If $(\vec{\phi},\,\Lambda)$ is an eigenpair of \eqref{e_36},
then so are $(\sigma_1\vec{\phi},\,-\Lambda)$, \ $(\vec{\phi}^*,\,-\Lambda^*)$, and 
$(\sigma_1\vec{\phi}^*,\,\Lambda^*)$, where
$$
\sigma_1 = \left( \ba{cc} 0 & 1 \\ 1 & 0 \ea \right)
$$
is another Pauli matrix. Note also that $\lambda$ is defined in the same way as in
\eqref{e_09}; hence Re$(\Lambda)\neq 0 $ indicates an instability.
Below we will use shorthand notations $\Lambda_R={\rm Re}(\Lambda)$ and 
$\Lambda_I={\rm Im}(\Lambda)$.

We begin analysis of \eqref{e_36} with two remarks. First, this equation is qualitatively
different from an analogous equation that arises in studies of stability of both bright
\cite{Kaup90} and dark \cite{Tran92} NLS solitons in that the relative sign of the first
and third terms of \eqref{e_36} is opposite of that in \cite{Kaup90,Tran92}. This fact
is the main reason why the unstable modes supported by \eqref{e_36} are qualitatively 
different from unstable modes of linearized NLS-type equations, as we will see below. 
While the latter modes are
supported by the soliton's core (see, e.g., Fig.~3 in \cite{Peli98}), the unstable
modes of \eqref{e_36} are supported by the soliton's ``tails". 

Second, from \eqref{e_36} and \eqref{e_35} one can easily establish the minimum value
of parameter $C$ where an instability (i.e., $\Lambda_R\neq 0 $) {\em can} occur.
The matrix operators on both sides of \eqref{e_36} are Hermitian; the operator 
$\sigma_3$ on the r.h.s. is not sign definite. Then the eigenvalues $\Lambda$
are guaranteed to be purely imaginary when the operator on the left-hand side (l.h.s.)
is sign definite \cite{RefLA}; otherwise they may be complex. The third term on the
l.h.s. of \eqref{e_36} is negative definite, and so is the first term in view of 
\eqref{e_34}. The second term, $D$, is negative when 
\be
C \, <  \, 1 /(|\beta| A^2).
\label{e_37}
\ee
Thus, \eqref{e_16} and \eqref{e_37} yield the stability condition of the fd-SSM on
the background of a soliton. We will show later that an unstable mode indeed
first arises when $C$ just slightly exceeds the r.h.s. of \eqref{e_37}.

Since the potential term in \eqref{e_36} is a slow function of $X$, it may seem
natural to employ the Wentzel--Kramers-Brillouin (WKB) method to analyze it.
Below we show that, unfortunately, the WKB method fails to yield an analytic form
of unstable modes of \eqref{e_36}. However, it still indicates {\em where} they can exist.
Away from ``turning points" (see below) the WKB-type solution of \eqref{e_36} is:
\be
\vec{\phi} = \left( a_+ e^{\theta_+/\epsilon} + b_+ e^{-\theta_+/\epsilon}\right)
 \, \vec{\varphi}_+ + 
 \left( a_- e^{\theta_-/\epsilon} + b_- e^{-\theta_-/\epsilon}\right)
 \, \vec{\varphi}_- ,
\label{e_38}
\ee
where $a_{\pm},\,b_{\pm}$ are some constants, and 
\bsube
\be
( \theta'_{\pm} )^2 = -D + 2V \pm \sqrt{V^2-\Lambda^2}, \qquad 
V \equiv V(\epsilon X), \;\; \theta' \equiv d\theta/d(\epsilon X),
\label{e_39a}
\ee
\be
\vec{\varphi}_{\pm} = 
 \frac{1}{\left[ ( \theta'_{\pm} )^2 (V^2-\Lambda^2) \right]^{1/4} } \,
 \left( \ba{r} \sqrt{\Lambda \pm \sqrt{\Lambda^2-V^2}} \vspace{0.2cm} \\
     -i\sqrt{\Lambda \mp \sqrt{\Lambda^2-V^2}}  \ea \right).
\label{e_39b}
\ee
\label{e_39}
\esube
At a turning point, say, $X=X_0$, the solution \eqref{e_38}, \eqref{e_39} breaks
down, which can occur because the denominator in \eqref{e_39b} vanishes.
In such a case, one needs to obtain a solution of \eqref{e_36} in a transition region
around the turning point by expanding the potential: \
$V(\epsilon X)=V(\epsilon X_0) + \epsilon(X-X_0)\,V'(\epsilon X_0)+\ldots$,
and then solving the resulting approximate equation. For a single Schr\"odinger
equation, a well-known solution of this type is given by the Airy function. This
solution is used to ``connect" the so far arbitrary constants $a_{\pm},\,b_{\pm}$
in \eqref{e_38} on both sides of the turning point.

Now a turning point of \eqref{e_36} is where: either (i) $\theta_+'=0$ or 
$\theta_-'=0$, \ or (ii) \ $(V(\epsilon X))^2-\Lambda^2 = 0$. The former case can be
shown (see, e.g., \cite{WKB_system}) to reduce to a single Schr\"odinger equation
case, where the solution in the transition region is given by the Airy function.
However, at present, no such transitional solution is analytically available in
case (ii)\footnote{Note that in this case, $(V^2-\Lambda^2)^{1/4}\vec{\varphi}_+$
and $(V^2-\Lambda^2)^{1/4}\vec{\varphi}_-$ are linearly dependent.}
\cite{Fulling2,Skorupski}.
Therefore, the solutions \eqref{e_38} canot be ``connected" by an analytic formula
across such a turning point, and hence one cannot find the eigenvalues $\Lambda$
in this case.

Based on the past experience with unstable linear modes of nonlinear waves, 
it is reasonable to assume that that unstable solutions of \eqref{e_36} must be 
localized. This allows one to predict the location of these solutions, which will
then be verified numerically. Indeed, expressions \eqref{e_38} and \eqref{e_39a}
show that $(\theta_{\pm}')^2>0$ for $(2V-D)>0$ and $V>\Lambda$ (i.e., inside the
soliton's core). The solution \eqref{e_38} cannot exponentially grow from, say,
the left side of the soliton to the right. Indeed, on the scale of Eq.~\eqref{e_36},
the soliton is very wide and hence, if the solution \eqref{e_38} is of order one
on the left, it would have become exponentially large on the right side of the soliton. 
Therefore, inside the soliton, the solution \eqref{e_38} must be exponentially
decaying. From \eqref{e_39a} it follows that it is also possible for a solution 
with $\Lambda_R\neq 0$ to decay outside the soliton (i.e., where $V\approx 0$)
when $D>0$. Thus, one can expect that a localized mode of \eqref{e_36} is to 
be lumped somewhere at the soliton's side and decay in both directions away from
that location.

\begin{figure}[h]
\vspace{-0.8cm}
\hspace*{-0.5cm}
\mbox{ 
\begin{minipage}{4cm}
\rotatebox{0}{\resizebox{4cm}{5.1cm}{\includegraphics[0in,0.5in]
 [8in,10.5in]{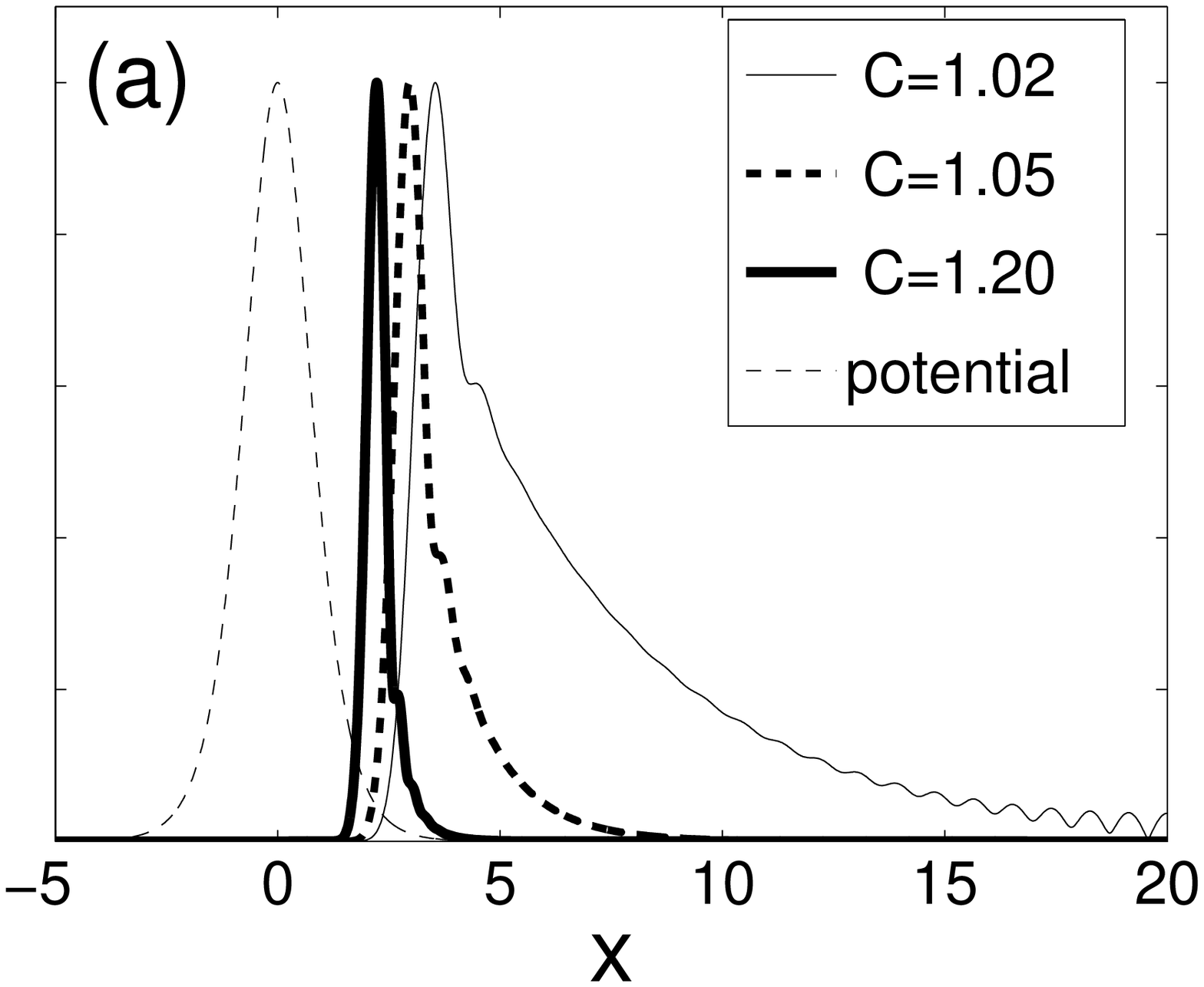}}}
\end{minipage}
\hspace{-0.4cm}
\begin{minipage}{4cm}
\rotatebox{0}{\resizebox{4cm}{5.1cm}{\includegraphics[0in,0.5in]
 [8in,10.5in]{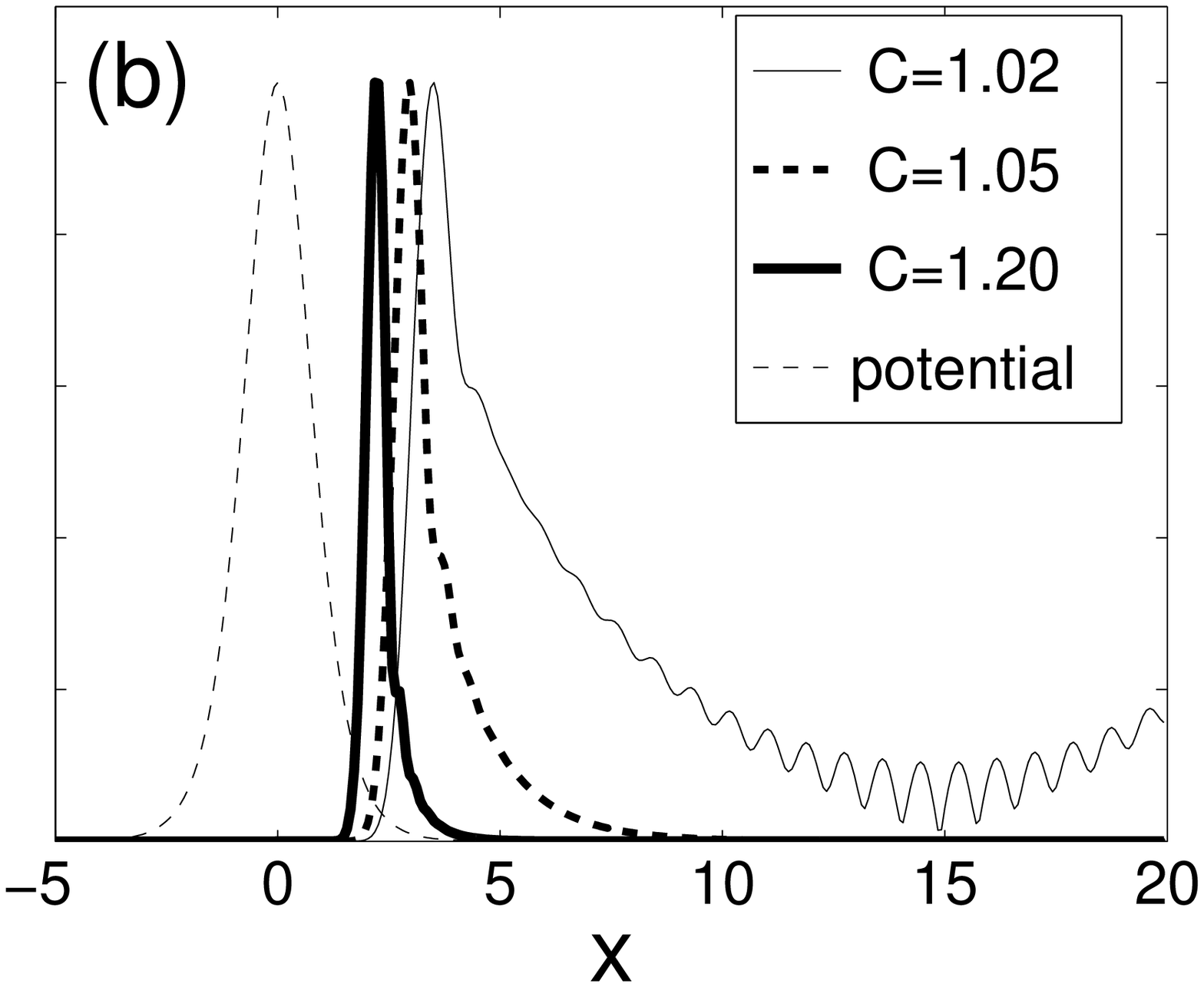}}}
\end{minipage}
\hspace{-0.4cm}
\begin{minipage}{4cm}
\rotatebox{0}{\resizebox{4cm}{5.1cm}{\includegraphics[0in,0.5in]
 [8in,10.5in]{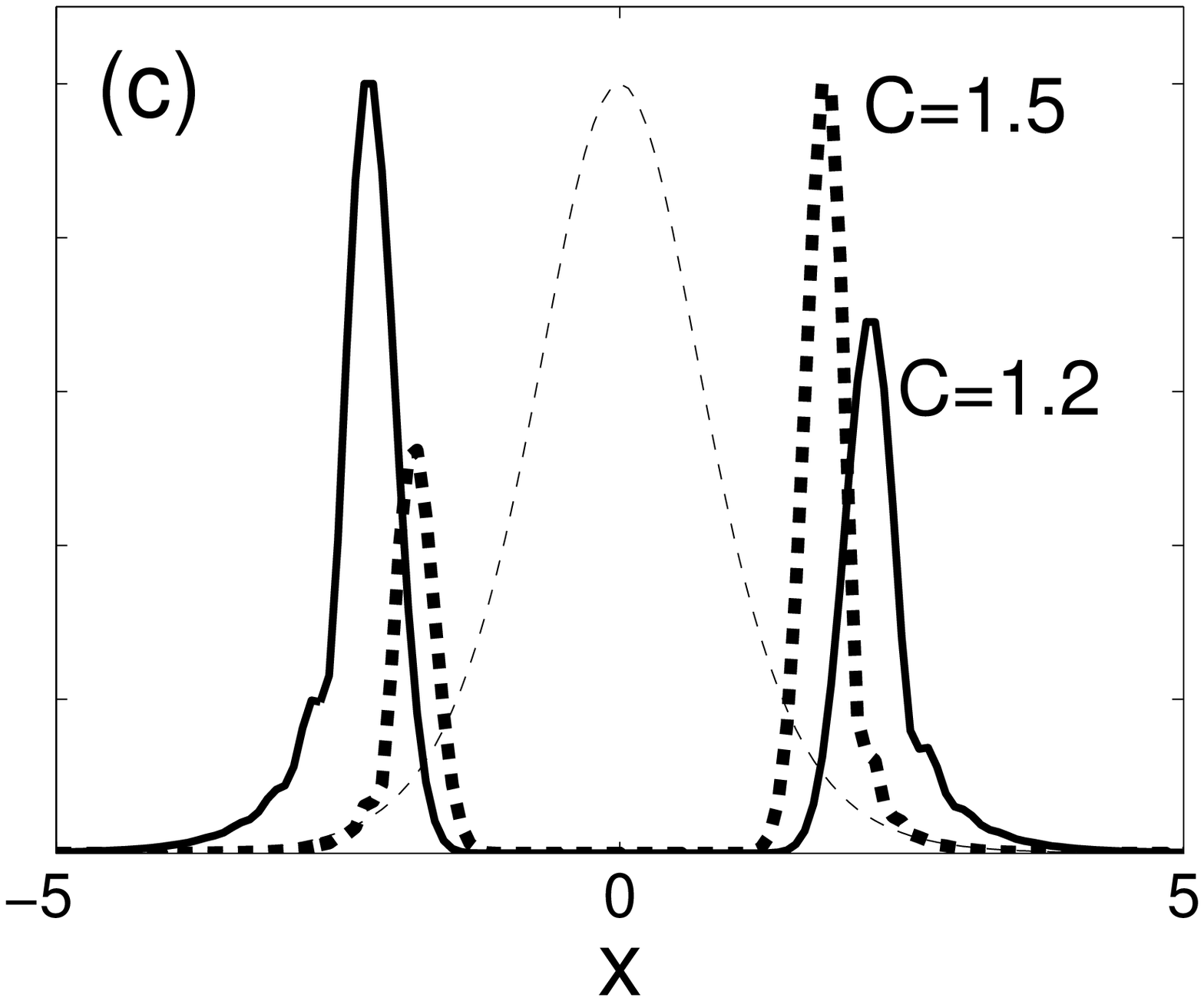}}}
\end{minipage}
\hspace{-0.3cm}
\begin{minipage}{3.8cm}
\rotatebox{0}{\resizebox{3.8cm}{5.1cm}{\includegraphics[0in,0.5in]
 [8in,10.5in]{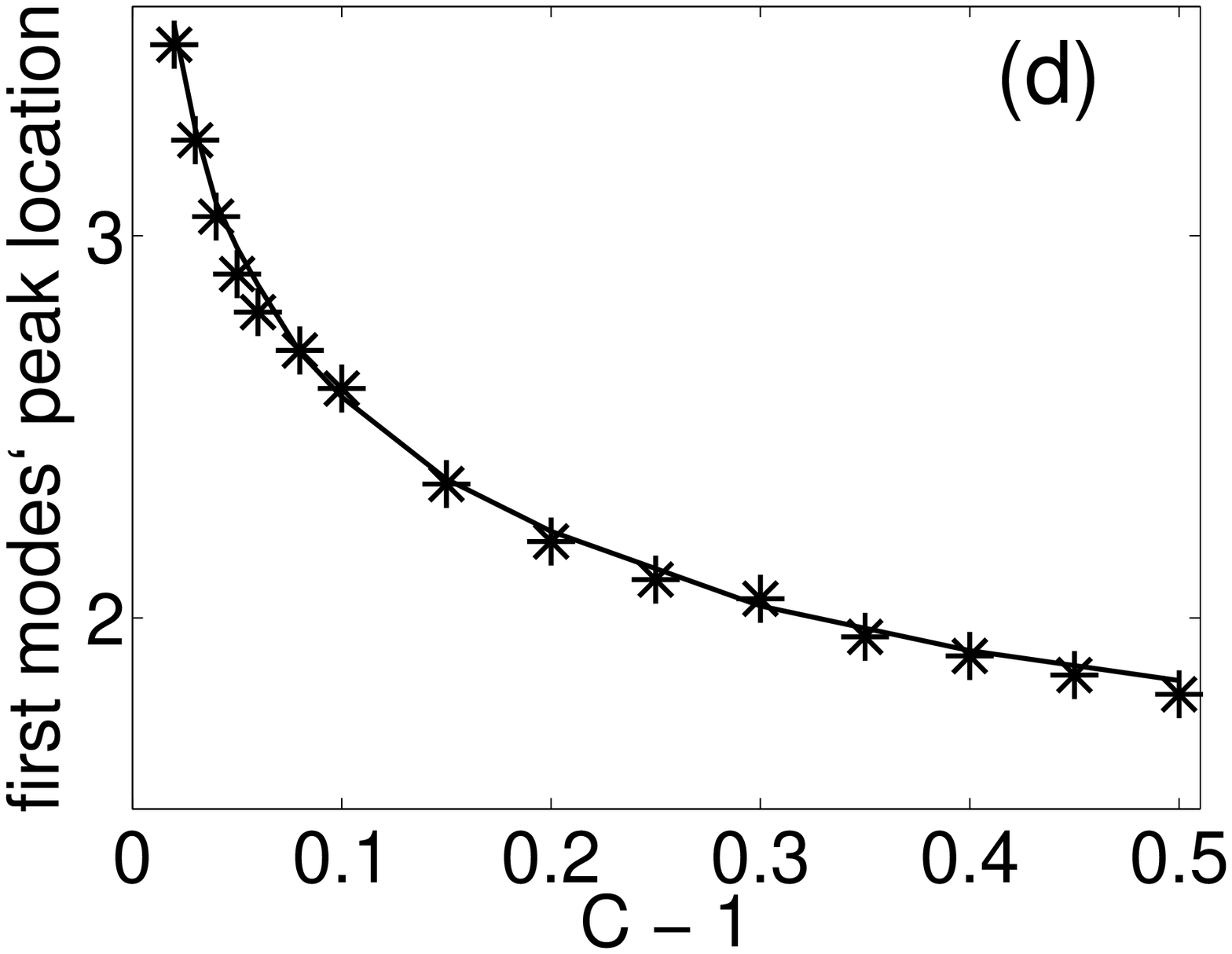}}}
\end{minipage}
 }
\vspace{-1.2cm}
\caption{(a): Profiles of the first localized mode on the right side of the soliton
for different values of $C$,
as found by the numerical method of Appendix B. \ (b): Same as in (a), but found 
from the numerical solution of \eqref{e_01}, as explained in the text. \ 
(c): The modes at {\em both} sides of the soliton found from the numerical solution
of \eqref{e_01}. Note that these modes do not ``see" each other because of the
barrier created by the soliton, and hence in general have different amplitudes as
they develop from independent noise seeds. In panels (a)--(c), the potential is 
$\sech^2(\epsilon X)$ (see \eqref{e_35}) and the amplitude of the mode is normalized
to that of the potential. \ (d): Location of the peak of the first localized mode,
found by the method of Appendix B (solid line) and from the solution of \eqref{e_01}
(stars). 
Similar data for $L=40$ and $N=2^{10}$ are very close to those in (d) and hence
are not shown.
}
\label{fig_5}
\end{figure}

In Fig.~\ref{fig_5} we show the first (i.e., corresponding to the greatest
$\Lambda_R$) such a mode for $L=40$, $N=2^9$ points (hence 
$\epsilon=\dx/2\approx 0.04$), $A=1$, $\beta=-1$, $\gamma=2$.
For these parameters, the threshold given by the r.h.s. of \eqref{e_37} is
$C=1$, and parameter $D$ in \eqref{e_36} is related to $C$ by:
\be
D=C-1.
\label{sec4_extra1}
\ee
The numerical method of finding these modes is described in Appendix B,
and the modes found by this method are shown in Fig.~\ref{fig_5}(a) for
different values of $C$. 
In Fig.~\ref{fig_5}(b) we show the same modes obtained from the numerical solution
of the NLS \eqref{e_01} by the fd-SSM. These
modes were extracted from the numerical solution by a high-pass filter, and then
the highest-frequency harmonic was factored out as per \eqref{e_26}. 
The agreement between Figs.~\ref{fig_5}(a) and \ref{fig_5}(b) is seen to be good.
In Fig.~\ref{fig_5}(d) we show the location of the peak of the first unstable mode,
computed both from \eqref{e_36} and from the numerical solution of \eqref{e_01},
versus parameter $C$. The corresponding values of the instability growth rate 
$\lambda$ were shown earlier in Fig.~\ref{fig_3}. Let us stress that $\lambda$
for the localized modes of \eqref{e_36} was found to be purely real 
up to the computer's round-off error ($\sim10^{-15}$).
There also exist unstable modes with complex $\lambda$, but such modes are not
localized and have smaller growth rates than the localized modes. 

\begin{figure}[h]
\vspace{-0.6cm}
\hspace*{-0.5cm}
\mbox{ 
\begin{minipage}{5.1cm}
\rotatebox{0}{\resizebox{5.1cm}{6.5cm}{\includegraphics[0in,0.5in]
 [8in,10.5in]{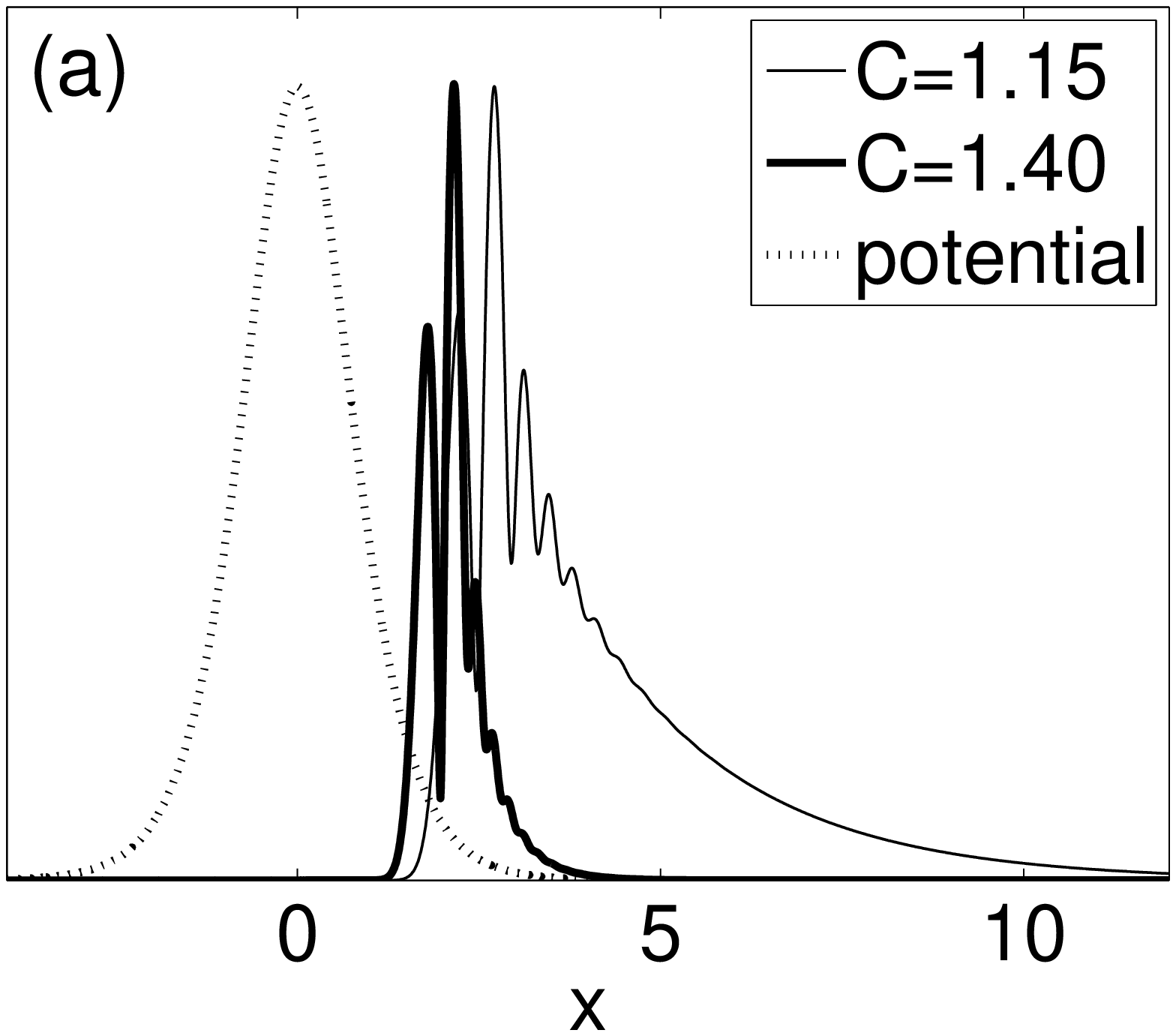}}}
\end{minipage}
\hspace{0.1cm}
\begin{minipage}{5.1cm}
\rotatebox{0}{\resizebox{5.1cm}{6.5cm}{\includegraphics[0in,0.5in]
 [8in,10.5in]{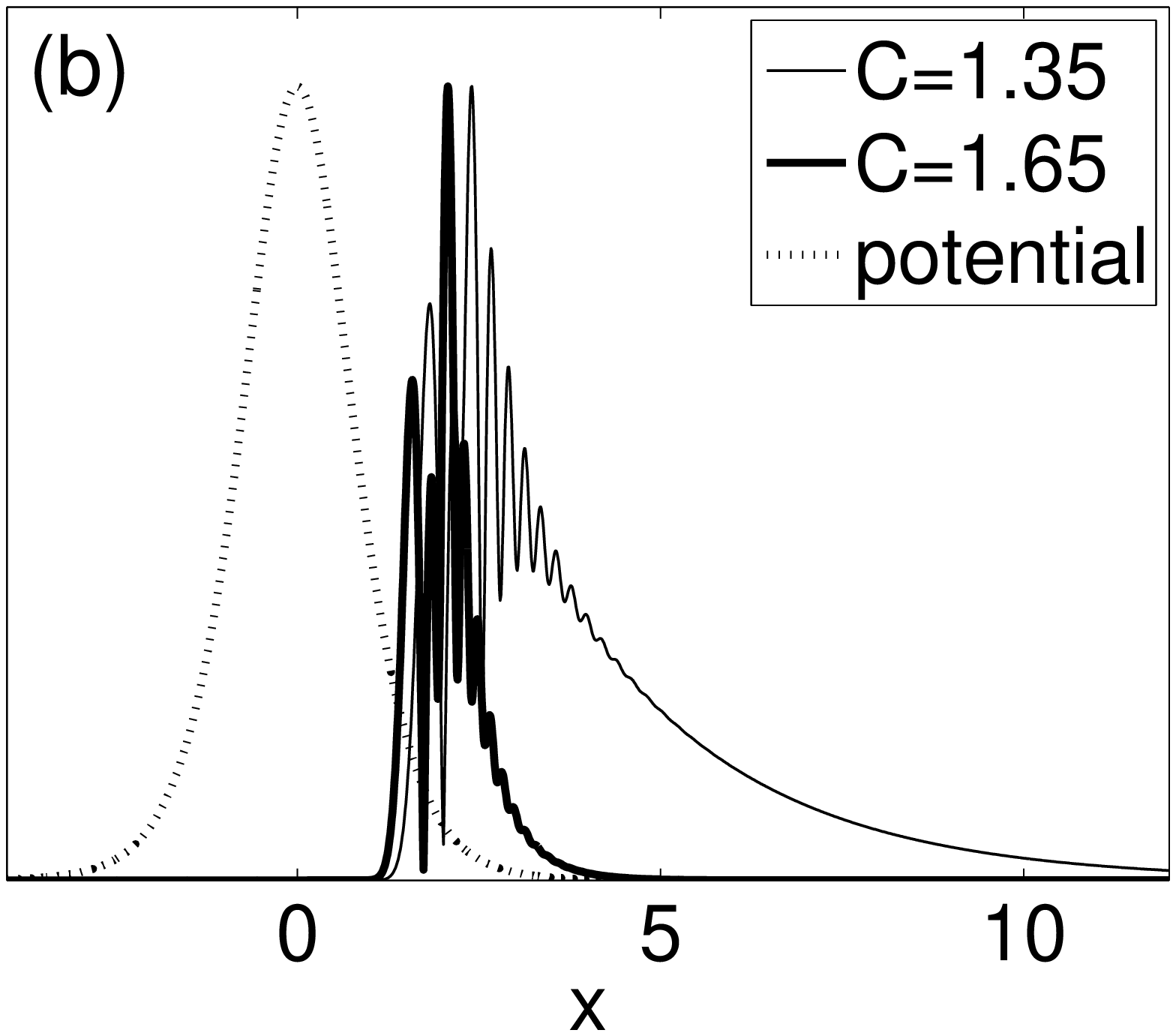}}}
\end{minipage}
\hspace{0.1cm}
\begin{minipage}{5.1cm}
\rotatebox{0}{\resizebox{5.1cm}{6.5cm}{\includegraphics[0in,0.5in]
 [8in,10.5in]{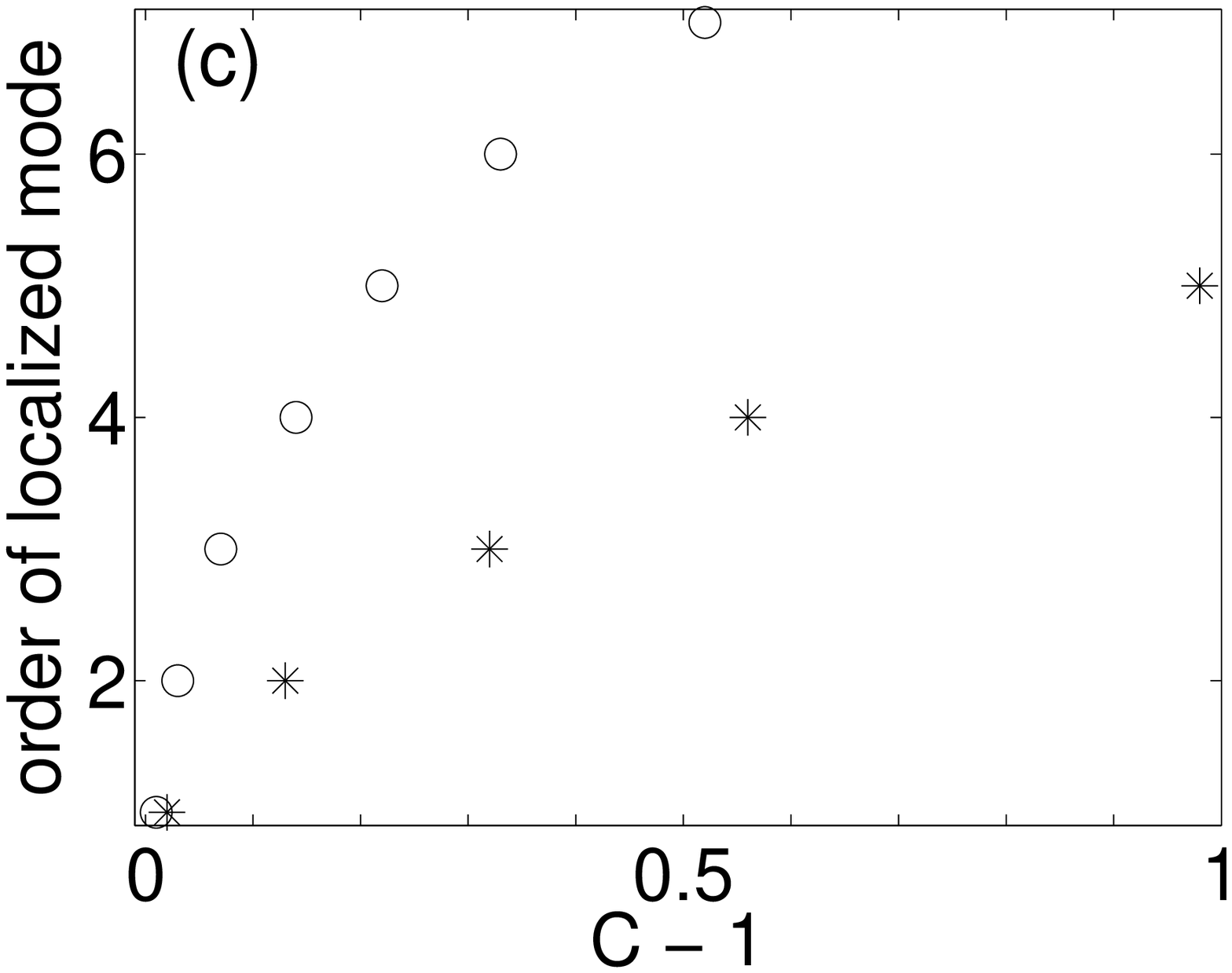}}}
\end{minipage}
 }
\vspace{-1.6cm}
\caption{Similar to Fig.~\ref{fig_5}(a), but for the second (a) and third (b) localized modes.
 \ (c): $C$ values where localized modes of increasing order appear. Stars --- for 
 $\epsilon = 40/1024$, circles --- for $\epsilon = 40/2048$. 
}
\label{fig_6}
\end{figure}

As $C$ increases from the critical value given by \eqref{e_37}, the localized unstable 
mode becomes narrower  and also moves toward the center of
the soliton. Moreover, higher-order localized modes of \eqref{e_36} arise. Typical
profiles of the second and third modes are shown in Fig.~\ref{fig_6}, along with the
parameter $C$ 
for which such modes first become localized within the spatial domain. 
In Appendix C we demonstrate that the process of ``birth" of an eigenmode that eventually
(i.e., with the increase of $C$) becomes localized,  is rather complicated. In particular,
it is difficult to pinpoint the exact value of parameter $C$ where such a mode appears.
Therefore, the $C$ values shown in Fig.~\ref{fig_6} are accurate only up to the second
decimal place.

\section{Conclusions}

The main contribution of this work is the extension the (in)stability analysis of the
fd-SSM beyond the von Neumann analysis in order to include spatially-varying background
solutions, such as the soliton \eqref{e_12}. We showed that, as previously for the
s-SSM \cite{ja}, this is done via a modified equation, \eqref{e_32}, derived for the
Fourier modes that approximately satisfy the resonance condition
\be
|\beta| k^2 \dt = \pi.
\label{conc_1}
\ee
Note that such modified equations --- \eqref{e_13} for the s-SSM and \eqref{e_32} for the fd-SSM ---
are different. Their analyses are also qualitatively different, and so are the modes that are
found to cause the instability of these two numerical methods. For the s-SSM, these modes are
almost monochromatic (i.e., non-localized) waves $\sim \exp(\pm ikx)$ that ``pass" through the
soliton very quickly. It is this scattering of those waves on the soliton that was shown \cite{ja} to 
lead to their instability. In contrast, for the fd-SSM considered in this work, the dominant unstable
modes are stationary relative to the soliton. Moreover, they are localized {\em at the sides}, as 
opposed to the core, of the soliton. To our knowledge, such localized modes were not reported before
in studies of (in)stability of nonlinear waves. We consider this finding another significant contribution
of this work.

Let us stress that, as previously for the s-SSM, the principle of frozen coefficients fails to even
qualitatively describe the instability of the fd-SSM on the soliton background. Indeed, that principle
assumes that the spatially  varying coefficient $U^2_{\rm b}(\epsilon \chi)$ 
in \eqref{e_32} can be approximated by
a constant. However, such an approach predicts \cite{WH} that any solution of the NLS with $\beta<0$,
including the soliton, should be stable (see 
\eqref{e_11} and the end of Sec.~III), contrary to the findings of this work.
In fact, the unstable modes of \eqref{e_36} localized at the soliton's sides owe their existence to the
spatial variation of the soliton's profile.

It was straightforward to obtain an approximate threshold, \eqref{e_37}
(where $C$ is given by \eqref{e_16}), beyond which the numerical
instability {\em may} occur. Our simulations showed that the instability does indeed occur just slightly
above that threshold. However, finding the exact threshold, $C_{\rm cr}$, remains an open problem, and so
is the calculation of the growth rate of the localized unstable modes; see Appendix C.

Let us note that the instability of the fd-SSM was studied in a recent paper \cite{Faou_2011}.
However, the focus of that paper was on a rigorous mathematical proof 
--- by a completely different method than we used here --- 
of the {\em stability} of the 
numerical scheme 
under a certain condition, rather than on the details of the development of the instability, as
in our work. That condition --- their Eq.~(2.9) --- is a constraint on the parameter
$r=\dt/\dx\,^2$, whereas our condition \eqref{e_37} is a constraint on the parameter $\dt/\dx$. 
Clearly, as $\dx\To 0$, our constraint allows for a greater $\dt$ than that of \cite{Faou_2011},
and hence is sharper.
Also, it is consistent with the numerical experiments reported in that paper.

Finally, let us show how our results can qualitatively explain the observed dynamics of the numerically
unstable soliton --- see the text after Eq.~\eqref{e_16}. Let $u_{\rm unst}$ be the field of the unstable
modes on the soliton's sides. 
At an early stage of the instabilty, $|u_{\rm unst}|\ll A$, the amplitude of the soliton. Also, its
characteristic wavenumbers are much greater than those of the soliton: see Fig.~\ref{fig_2}(a)
and \eqref{e_23}.
Accordingly, substitute $u=u_{\rm sol}+u_{\rm unst}$ into the NLS \eqref{e_01} and discard all the 
high-wavenumber terms to obtain:
\be
i(\usol)_t -\beta (\usol)_{xx}+\gamma \usol |\usol|^2 = - 2\gamma \usol |u_{\rm unst}|^2.
\label{conc_2}
\ee
This is the equation for a perturbed soliton with the perturbation being, in general, not symmetric about
the soliton's center (see Fig.~\ref{fig_5}(c)), because the modes on the left and right sides 
of the soliton do not ``see" each other and hence can have different amplitudes. Such a perturbation is
known (see, e.g., \cite{Agrawal_book}, Sec.~5.4.1) to cause the soliton to move, which is precisely the 
effect we observed.


\section*{Acknowledgement}

I thank Jake Williams for help with numerical simulations at an early stage of this work,
and Eduard Kirr for a useful discussion.
This research was supported in part by the NSF grant ECCS-0925706.


\section*{Appendix A: 
Modified linearized NLS for fd-SSM with non-periodic boundary conditions}

We consider the Dirichlet boundary conditions (b.c.)
and without loss of generality (for the
purposes of this derivation) set them equal to zero. Neumann or mixed b.c.
can be treated similarly, and lead to similar results. 

The equation for the dispersive step of the SSM is still given by \eqref{e_05}.
However, now instead of \eqref{e_06} we assume: $u_{n+1}^0=0$, $u_{n+1}M=0$, where
$m=0$ and $m=M$ are the end points of the spatial grid. Then \eqref{e_05} can be
rewritten as \cite{KincaidCheney}
\be
({\mathcal I}+(i\beta r/2){\mathcal A})\, {\bf u}_{n+1} = 
({\mathcal I}-(i\beta r/2){\mathcal A})\, {\bf \bu},
\label{e_40}
\ee
where: ${\bf \bu}=[\bu^1,\bu^2,\ldots,\bu^{M-1}]^T$, similarly for ${\bf u}_{n+1}$,
${\mathcal I}$ is an $(M-1)\times (M-1)$ identity matrix, and ${\mathcal A}$ is an
$(M-1)\times (M-1)$ tridiagonal matrix with $(-2)$ on the main diagonal and $(+1)$
on the sub- and super-diagonals. 

The starting point of our derivation in Sec.~III, Eq.~\eqref{e_18}, has exacly the 
same form for the case of the Dirichlet b.c., except that $\F$
is replaced with ${\mathcal T}$ --- an expansion over the complete set of the
eigenvectors of ${\mathcal A}$; similarly, $\F^{-1}$ is replaced by ${\mathcal T}^{-1}$.
The exponential in \eqref{e_19} that acts on the $j$th eigenvector is replaced by
\be
e^{iP_j}=\frac{1-i\beta r\lambda_j/2}{1+i\beta r\lambda_j/2},
\label{e_41}
\ee
where $\lambda_j$ is the corresponding eigenvalue \cite{KincaidCheney}:
\be
\lambda_j = -4 \sin^2\big( \pi j/(2M) \big).
\label{e_42}
\ee
Equations \eqref{e_41}, \eqref{e_42} and the middle expression in \eqref{e_19}
coincide provided that we identify:
\be
k=j\pi/(M\dx) = j\pi/L.
\label{e_43}
\ee
However, we are still a step away from proving that the modified linearized NLS for
the Dirichlet b.c. case is the same as that equation for periodic b.c..
This is because $-k^2$, which is the Fourier symbol of the second derivative, is
not the symbol of the second derivative under the transformations ${\mathcal T}$
and ${\mathcal T}^{-1}$. Under those transformation, the required symbol is given 
by \eqref{e_42}. We will now use this observation to supply the last step and show
that the modified linearized NLS for the case of Dirichlet b.c. is indeed the same
as \eqref{e_31}. This follows from \eqref{e_41}--\eqref{e_43} and a calculation that
is similar to \eqref{e_21}:
\bea
e^{iP(k)} & \approx & -\left( 1 + \frac1{i\beta r 
                                         \big(1-\sin^2((k-k_{\max})\dx/2)\big)} \right)
 \nonumber \\
& \approx & -\left( 1 + \frac1{i\beta r} +
 \frac{\sin^2((k-k_{\max})\dx/2)\big)}{i\beta r} \right),
\label{e_44}
\eea
where we have used that $\sin(k_{\max}\dx/2)=1$ and that for highly oscillatory
eigenvectors of ${\mathcal A}$, one has $(k-k_{\max})\dx \ll 1$. The last term
on the r.h.s. of \eqref{e_44} is the desired symbol of the second derivative, and 
then the rest of the derivation is the same as that leading to \eqref{e_28}. From it
one obtains the same modified linearized NLS as \eqref{e_31}. Our numerical simulations
of the NLS using the fd-SSM with zero Dirichlet b.c. confirm this conclusion.


\section*{Appendix B: 
Numerical solution of eigenproblem \eqref{e_36}}

We work with \eqref{e_36} written in an equivalent form:
\be
\sigma_3 \left( \partial_X^2 + D - V(\epsilon X) 
 \left( \ba{cc} 2 & 1 \\ 1 & 2 \ea \right) -i\Lambda_0 \sigma_3 \,\right)
\vec{\phi} \,=\, i(\Lambda-\Lambda_0) \sigma_3 \vec{\phi},
\label{B_01}
\ee
where the reason to include a constant $\Lambda_0$ will be explained later. 
We discretize \eqref{B_01} using Numerov's method. It approximates the equation \ 
$\Phi_{XX}=F(\Phi,X)$ \ by a finite-difference scheme
\be
\Phi^{m+1}-2\Phi^m+\Phi^{m-1}=\frac{\D X\,^2}{12} \big( F^{m+1}+10F^m+F^{m-1} \big)
\label{B_02}
\ee
with accuracy $O(\D X \,^4)$; here $\Phi^m\equiv \Phi(X_m)$, $F^m\equiv F(\Phi^m,X_m)$, etc., and
$m=1,\ldots\,,M-1$. Then, for the discretized solution ${\bf f}_k=[\phi_k^1,\ldots\,,\,\phi_k^{M-1}]^T$
($k=1,2$) one obtains:
\be
(-1)^{k-1}\left(\, \left[ \frac1{\D X\,^2}{\mathcal A}_{\rm per} + {\mathcal N}_{\rm per} 
\left\{ D{\mathcal I} - 2{\mathcal V} - i\Lambda_0 (-1)^{k-1} {\mathcal I} \right\} \right] {\bf f}_k
\,-\, {\mathcal N}_{\rm per} {\mathcal V} {\bf f}_{3-k} \, \right) \,=\, 
i(\Lambda - \Lambda_0) {\mathcal N}_{\rm per} {\bf f}_k,
\label{B_03}
\ee
where ${\mathcal I}$ is defined after \eqref{e_40}; \ ${\mathcal A}_{\rm per}$ is as in \eqref{e_40}
except that its $(1,M-1)$th and $(M-1,1)$th entries equal 1 (to account for the periodic boundary
conditions); \ ${\mathcal N}_{\rm per}$ has a similar structure as ${\mathcal A}_{\rm per}$: \ 
$({\mathcal N}_{\rm per})_{m,m}=10/12$ (see \eqref{B_02}), $({\mathcal N}_{\rm per})_{(m-1),m}=
({\mathcal N}_{\rm per})_{m,(m-1)}=1/12$, $({\mathcal N}_{\rm per})_{1,(M-1)}=
({\mathcal N}_{\rm per})_{(M-1),1}=1/12$, and the rest of its entries are zero; \ 
and ${\mathcal V}={\rm diag}(V^1,\ldots\,,\,V^{M-1})$. 
Next, defining the combined vector and matrices:
$$
\hat{\bf f}=\left[ \ba{c} {\bf f}_1 \\ {\bf f}_2 \ea \right], \quad
\hat{\mathcal A}_{\rm per}= \left( \ba{cc} {\mathcal A}_{\rm per} & {\mathcal O} \\
                                   {\mathcal O} & {\mathcal A}_{\rm per} \ea \right), \quad
\hat{\mathcal N}_{\rm per}= \left( \ba{cc} {\mathcal N}_{\rm per} & {\mathcal O} \\
                                   {\mathcal O} & {\mathcal N}_{\rm per} \ea \right), 
$$
$$
\hat{\mathcal V}= \left( \ba{cc} 2{\mathcal V} & {\mathcal V} \\
                                   {\mathcal V} & 2{\mathcal V} \ea \right), \quad
\hat{\sigma}_3 = \left( \ba{cc} {\mathcal I} & {\mathcal O} \\
                                   {\mathcal O} & -{\mathcal I} \ea \right),
$$
where ${\mathcal O}$ is the $(M-1)\times (M-1)$ zero matrix, one rewrites \eqref{B_03} as:
\be
\hat{\sigma}_3 
\left[ \frac1{\D X\,^2}\hat{\mathcal A}_{\rm per} + D\hat{\mathcal N}_{\rm per} 
 - \hat{\mathcal N}_{\rm per} \hat{\mathcal V} - i\Lambda_0 \hat{\sigma}_3 \hat{\mathcal N}_{\rm per} 
  \right] \hat{\bf f} \,=\, i(\Lambda-\Lambda_0)\hat{\mathcal N}_{\rm per} \hat{\bf f}\,.
\label{B_04}
\ee
This equation has the form of the generalized eigenvalue problem 
${\mathcal G}\hat{\bf f}=\lambda{\mathcal H}\hat{\bf f}$ where ${\mathcal H}=\hat{\mathcal N}_{\rm per}$
is a positive definite matrix. This problem can be solved by Matlab's command \verb+eigs+.
As its options, we specified that we were looking for 24 smallest-magnitude eigenvalues and
the corresponding eigenmodes. Beyond the instability threshold there are several modes with complex $\Lambda$.
We visually inspected them and found that the most unstable mode was also the most localized and also had 
a real eigenvalue.

We verified that the eigenvalues did not change to five significant figures whether we used $\D X=1/10$
or $1/20$; so we used $\D X=1/10$. Finally, the need to shift the eigenvalues by a $\Lambda_0$ occurred
for relatively large $D$: $D\ge 0.1$. There, the localized modes were no longer among those with the 24 
smallest-magnitude eigenvalues (see above), and significantly increasing the number of modes sought
caused \verb+eigs+ to fail. Therefore, we had to shift the eigenvalues by $\Lambda_0$ in order to make
the mode of interest appear among the 24 smallest-eigenvalue ones.


\section*{Appendix C: 
``Birth" of localized unstable mode}

The instability growth rates plotted in Fig.~\ref{fig_3} are monotonic functions of the
the parameter $C$. This, however, occurs only when $C$ is sufficiently beyond a critical
value, $C_{\rm cr}$, where the dominant (i.e., with the greatest-$|\Lambda_R|$) unstable 
mode is created. Near $C_{\rm cr}$, which
is slightly above the threshold value given by the r.g.s. of \eqref{e_37}, the evolution
of the greatest-$|\Lambda_R|$ eigenvalue is quite irregular.

Below we present results about this evolution for the dominant unstable mode (shown in
Fig.~\ref{fig_5}) for $L=40$, $N=2^9$ (i.e., $\epsilon=40/1024\approx 0.04$) and the rest of
the parameters being the same as listed in Sec.~II, i.e.: $\beta=-1$, $\gamma=2$, and $A=1$. 
Then, the control parameter $C$ is related to the free term in \eqref{e_36} by \eqref{sec4_extra1}.
While we have been unable to rigorously establish an analytical expression for $C_{\rm cr}$,
we will present a hypothesis as to what it may be. The main message that we intend to convey
is that the ``birth" of a localized eigenmode of Eq.~\eqref{e_36} occurs via a complex 
sequence of bifurcations, in contrast to a single bifurcation that typically takes place
when an unstable mode of a nonlinear wave is ``born" (see, e.g., \cite{Kapitula98}).

We numerically  observed that eigenmodes of \eqref{e_36} with $\Lambda_R\neq 0$ emerge
not only from the origin (i.e. when imaginary eigenvalues ``collide" at $\Lambda=0$
and give rise to two real ones),
but also from ``collision" of eigenvalues with $\Lambda_I\neq 0$ whereby two\footnote{
As per the remark after \eqref{e_36}, there is also a pair of eigenvalues with $(-\Lambda_I)$,
so in this case there exists a quadruplet of complex eigenvalues.}
complex eigenvalues are ``born". However, the very first (i.e., for the smallest $C$)
unstable mode does emerge out of the origin. We will now show that this mode is essentially
non-localized and, moreover, it is {\em not} the mode that eventually becomes the dominant
unstable mode, whose growth rate is plotted in Fig.~\ref{fig_3}. The reason why we still chose
to discuss the former mode while being primarily interested in the latter one, will become clear
as we proceed.

\begin{figure}[h]
\vspace{-1.6cm}
\mbox{ 
\begin{minipage}{7cm}
\rotatebox{0}{\resizebox{7cm}{9cm}{\includegraphics[0in,0.5in]
 [8in,10.5in]{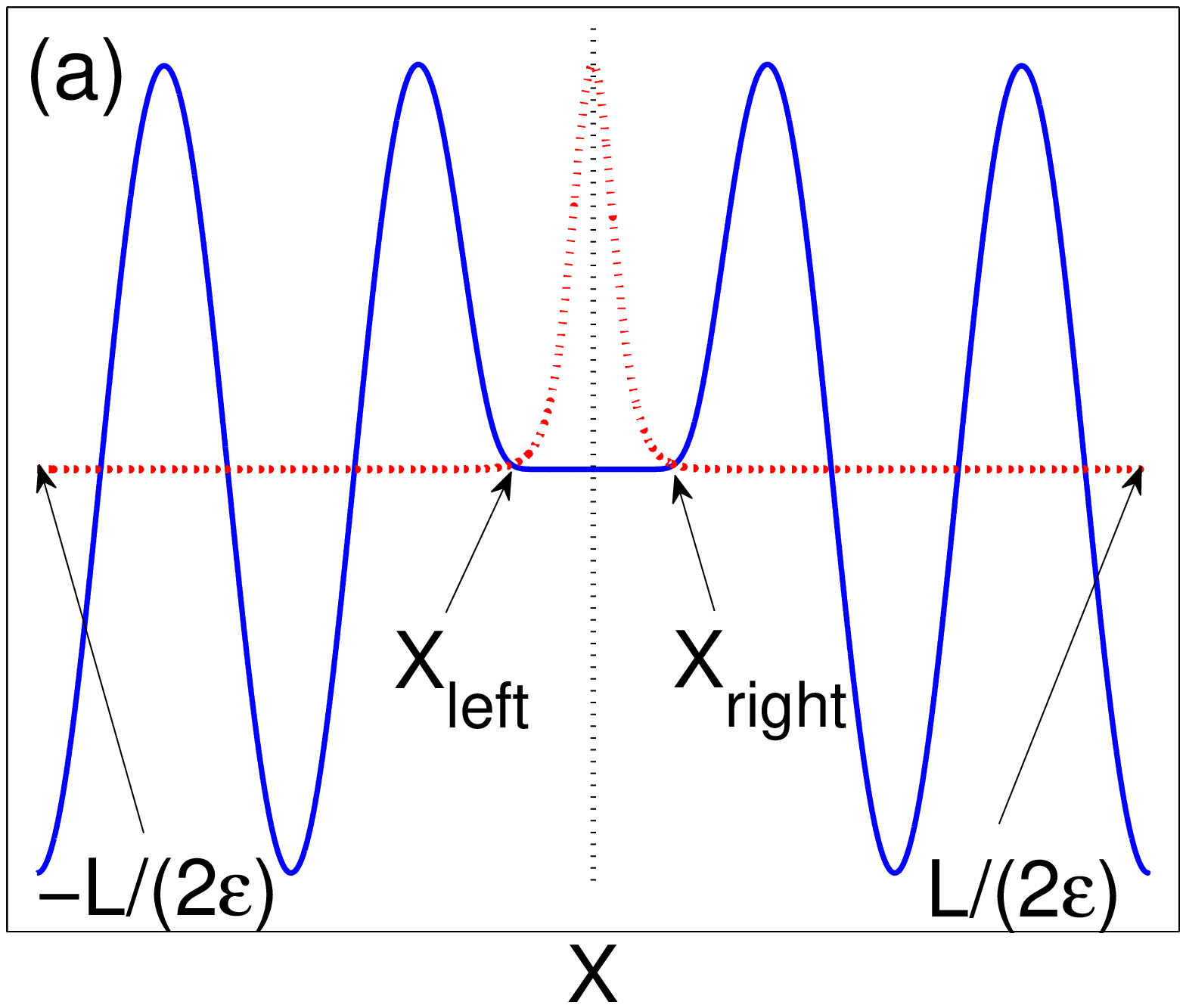}}}
\end{minipage}
 \hspace{0.1cm}
 \begin{minipage}{7cm}
  \rotatebox{0}{\resizebox{7cm}{9cm}{\includegraphics[0in,0.5in]
   [8in,10.5in]{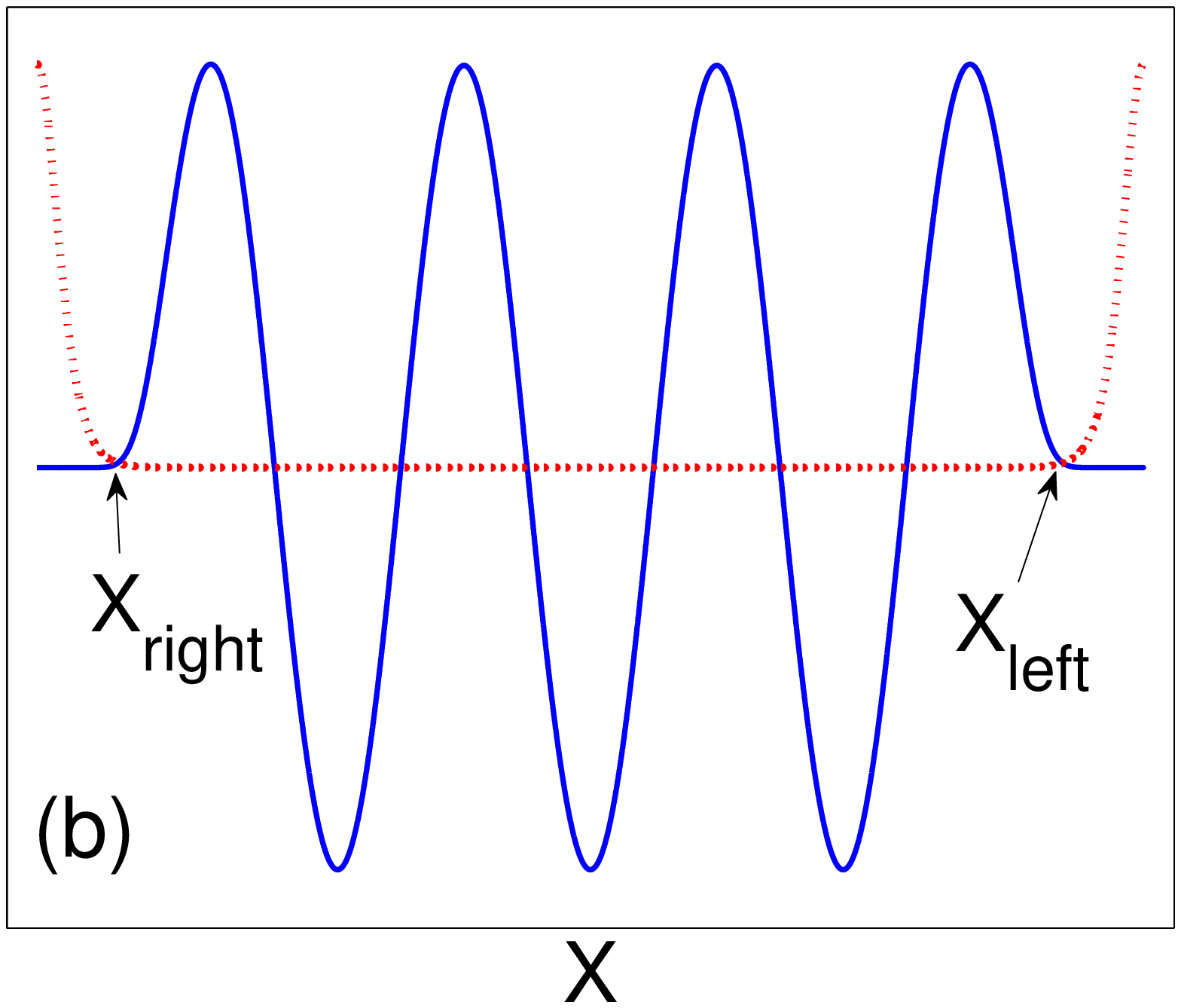}}}
 \end{minipage}
 }
\vspace{-1.6cm}
\caption{ (Color online) \ (a): A solution of \eqref{e_C1} (solid); potential $\sech^2(\epsilon X)$ 
(red dotted). The amplitude of the solution is normalized to that of the potential. \ 
(b): Same as (a), but that panel is ``cut" along the vertical dotted line at the center,
and the resulting halves are interchanged.
}
\label{fig_7}
\end{figure}

For a mode with $\Lambda=0$, Eq.~\eqref{e_36} can be split into two uncoupled Schr\"odinger
equations:
\bsube
\be
\big( \,\partial_X^2 + D - \nu_{\pm}\,V(\epsilon X)\,\big) \phi_{\pm} =0, \qquad 
\phi_{\pm} = \phi_1 \pm \phi_2,
\label{e_C1a}
\ee
where $\nu_-=1$ and $\nu_+=3$. Note that $\phi_{\pm}$ satisfy the periodic boundary conditions:
\be
\phi_{\pm}(-L/(2\epsilon)) \,=\, \phi_{\pm}(L/(2\epsilon)) .
\label{e_C1b}
\ee
\label{e_C1}
\esube
In Fig.~\ref{fig_7}(a) we show an example of a nontrivial solution of \eqref{e_C1}.
In view of the periodic boundary conditions, this figure is equivalent to Fig.~\ref{fig_7}(b).
Recall from Sec.~IV that the eigenmode is exponentially small inside the soliton. Then
the solution shown in Fig.~\ref{fig_7}(b) can be thought of as being localized inside the valley
bounded by the two ``halves" of the potential. Using this observation, one can estimate
the isolated values of $D$ for which one of the equations \eqref{e_C1a}, along with
\eqref{e_C1b}, has a nontrivial solution, by the WKB method. The condition for the existence
of a mode localized inside the valley of Fig.~\ref{fig_7}(b) is given by the Bohr--Sommerfeld
formula:
\be
\left( \int_{-L/(2\epsilon)}^{X_{\rm left}} + \int_{X_{\rm right}}^{L/(2\epsilon)} \right)
 \, \sqrt{D-\nu V(\epsilon X)} \,dX \,=\, \pi \left(n+\frac12 \right),
\label{e_C2}
\ee
where $\nu$ is either $\nu_-$ or $\nu_+$, $n$ is an integer, and $X_{\rm left,\,right}$ are
the turning points (see Sec.~IV), where 
\be
D-\nu V(\epsilon X_{\rm left,\,right}) =0.
\label{e_C3}
\ee
The number of full oscillation periods of the mode inside the valley equals $n$; for example,
in Fig.~\ref{fig_7}, $n=3$.

When $D\ll 1$, the $\sech$ potential in \eqref{e_C3} can be approximated by an exponential: \ 
$\sech^2(\epsilon X) \approx 4\exp(-2\epsilon X)$. Then, using \eqref{e_35} and \eqref{e_C3},
we reduce \eqref{e_C2} to
\be
\sqrt{D} 
\int_{X_{\rm right}}^{L/(2\epsilon)} \sqrt{1 - \exp[-2\epsilon(X-X_{\rm right})] } \,dX \,=\,
\frac{\pi}2 \left( n+\frac12 \right),
\label{e_C4}
\ee
with $X_{\rm left}=-X_{\rm right}$ and 
\be
X_{\rm right} = \frac1{2\epsilon} \, \ln \frac{8\nu C \beta^2}{D}\,.
\label{e_C5}
\ee
Neglecting the exponentially small terms of the order $O\big( \, \exp[-(L-2\epsilon X_{\rm right})]\,\big)$,
one obtains from \eqref{e_C4}:
\be
\sqrt{D} \left( L - \ln \frac{8\nu C\beta^2}{D} - 2(1-\ln 2) \right) = \epsilon\,\pi \left( n+\frac12 \right).
\label{e_C6}
\ee

Note that the WKB condition \eqref{e_C2}, and hence \eqref{e_C6}, is valid when $n$ is sufficiently
large. In particular, it is {\em not} supposed to accurately predict the ``birth" of the first 
unstable mode, where $n=0$. Indeed, Eq.~\eqref{e_C6} predicts that such a mode (for $\nu=1$)
emerges at $D\approx 5.9\cdot 10^{-6}$, while numerically (see Appendix B) it is found at 
$D\approx 1.6\cdot 10^{-5}$. (A similar mode for $\nu=3$ emerges at a slightly higher value
of $D$.) Formula \eqref{e_C6} becomes accurate to the fourth significant figure in $D$ for
$n \gtrsim 20$.

As we noted above, the first unstable mode is {\em not} the one that eventually becomes the 
dominant unstable mode. It disappears already at $D\approx 1.7\cdot 10^{-5}$, and there is
an adjacent interval of $D$ values where all the eigenvalues of \eqref{e_36} are purely imaginary
(i.e., the soliton is numerically stable). 
As $D$ increases, higher-order ``real" (i.e., with $\Lambda_I=0$) modes appear and disappear in a 
similar fashion, as do quadruplets of modes with complex $\Lambda$. In both these types of modes,
$\Lambda_R$ is fairly small: $|\Lambda_R| \lesssim D/10$. There also exist intervals of $D$,
of increasingly small length, where all $\Lambda$'s are purely imaginary. This situation persists
until the dominant unstable mode appears at $D_{\rm cr}\,(=C_{\rm cr}-1)$. This occurs as follows.

\begin{figure}[h]
\vspace*{-0.6cm}
\hspace*{-0.5cm}
\mbox{ 
\begin{minipage}{5.1cm}
\rotatebox{0}{\resizebox{5.1cm}{7cm}{\includegraphics[0in,0.5in]
 [8in,10.5in]{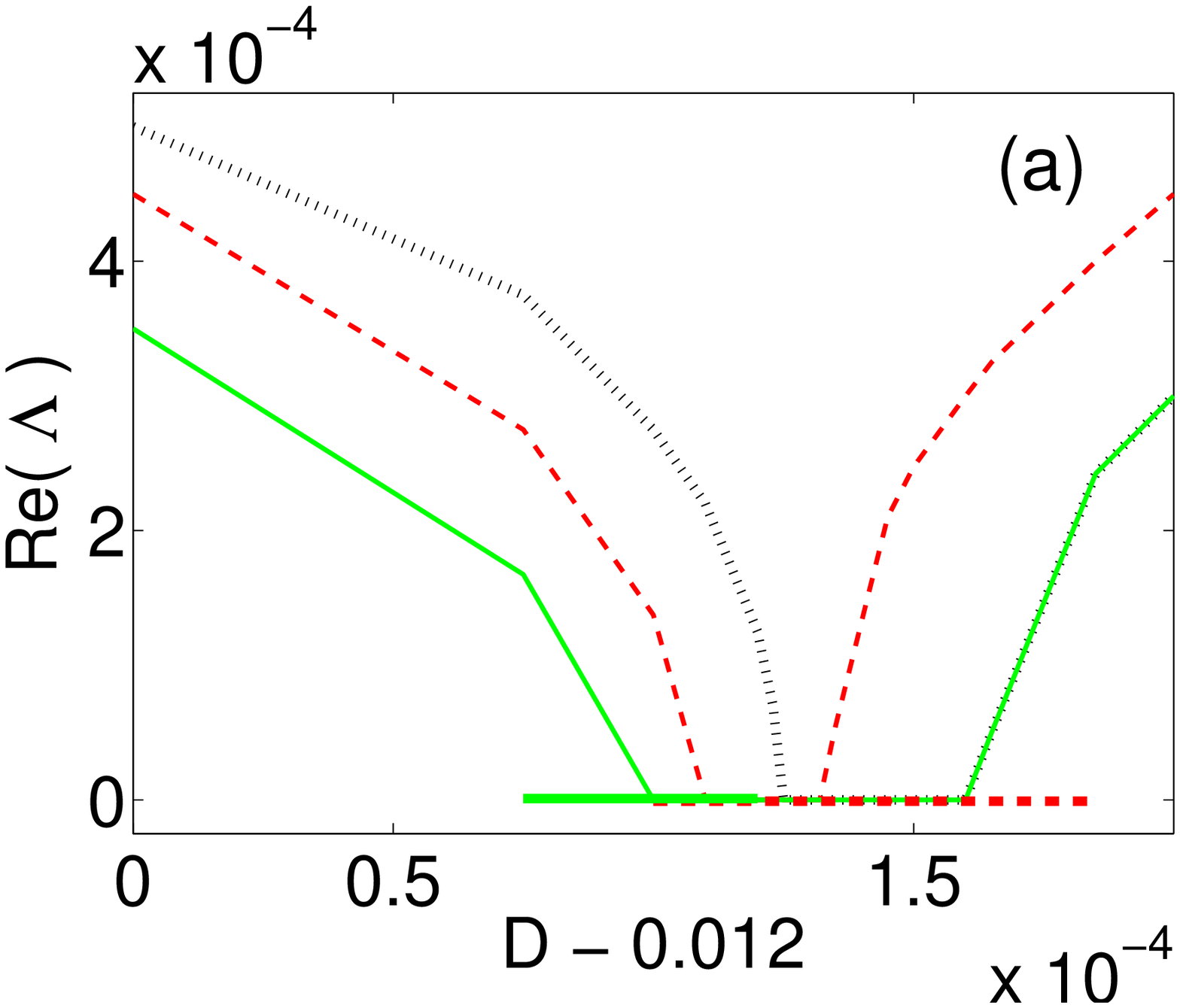}}}
\end{minipage}
\hspace{0.1cm}
\begin{minipage}{5.1cm}
\rotatebox{0}{\resizebox{5.1cm}{7cm}{\includegraphics[0in,0.5in]
 [8in,10.5in]{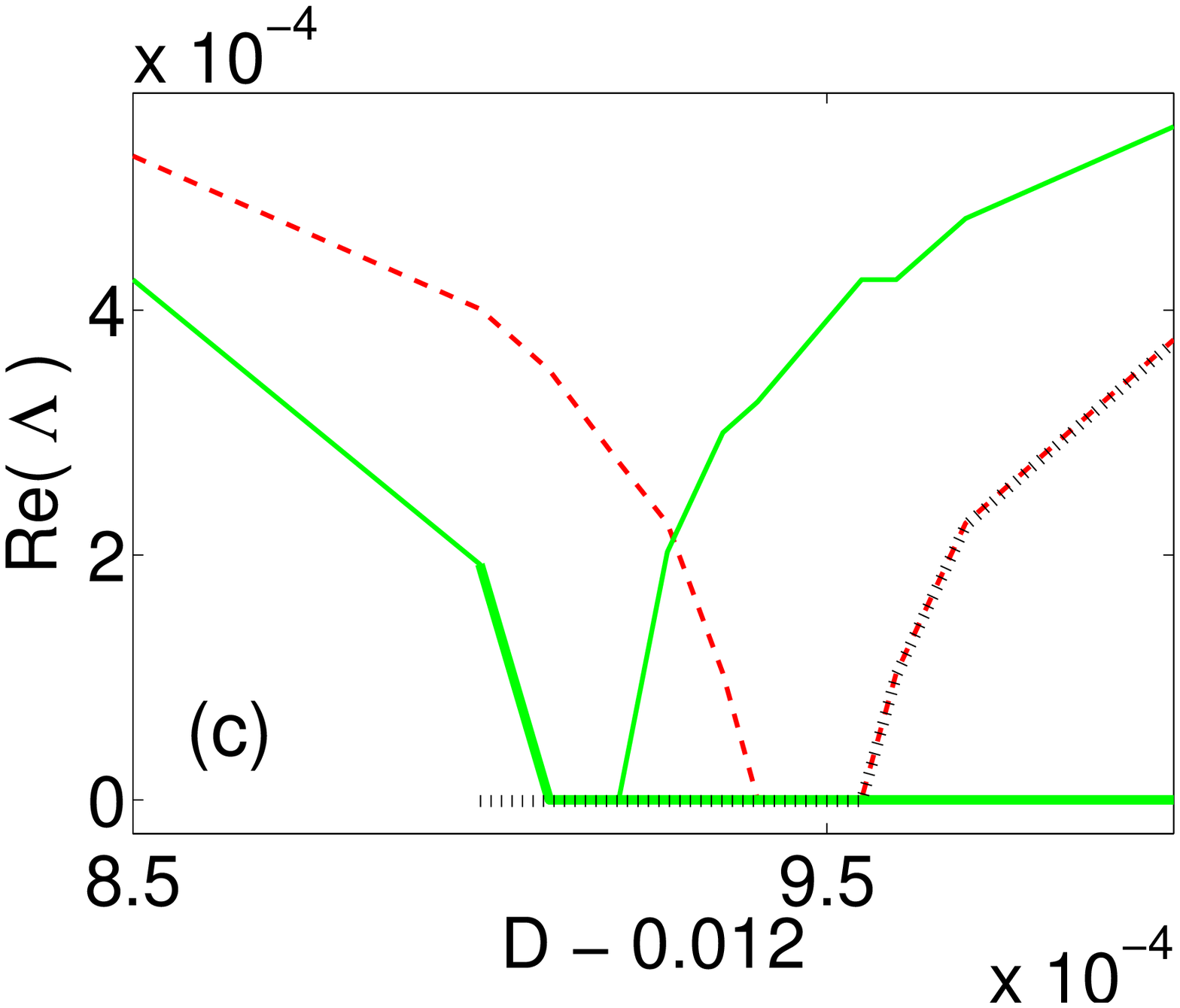}}}
\end{minipage}
\hspace{0.1cm}
\begin{minipage}{5.1cm}
\rotatebox{0}{\resizebox{5.1cm}{7cm}{\includegraphics[0in,0.5in]
 [8in,10.5in]{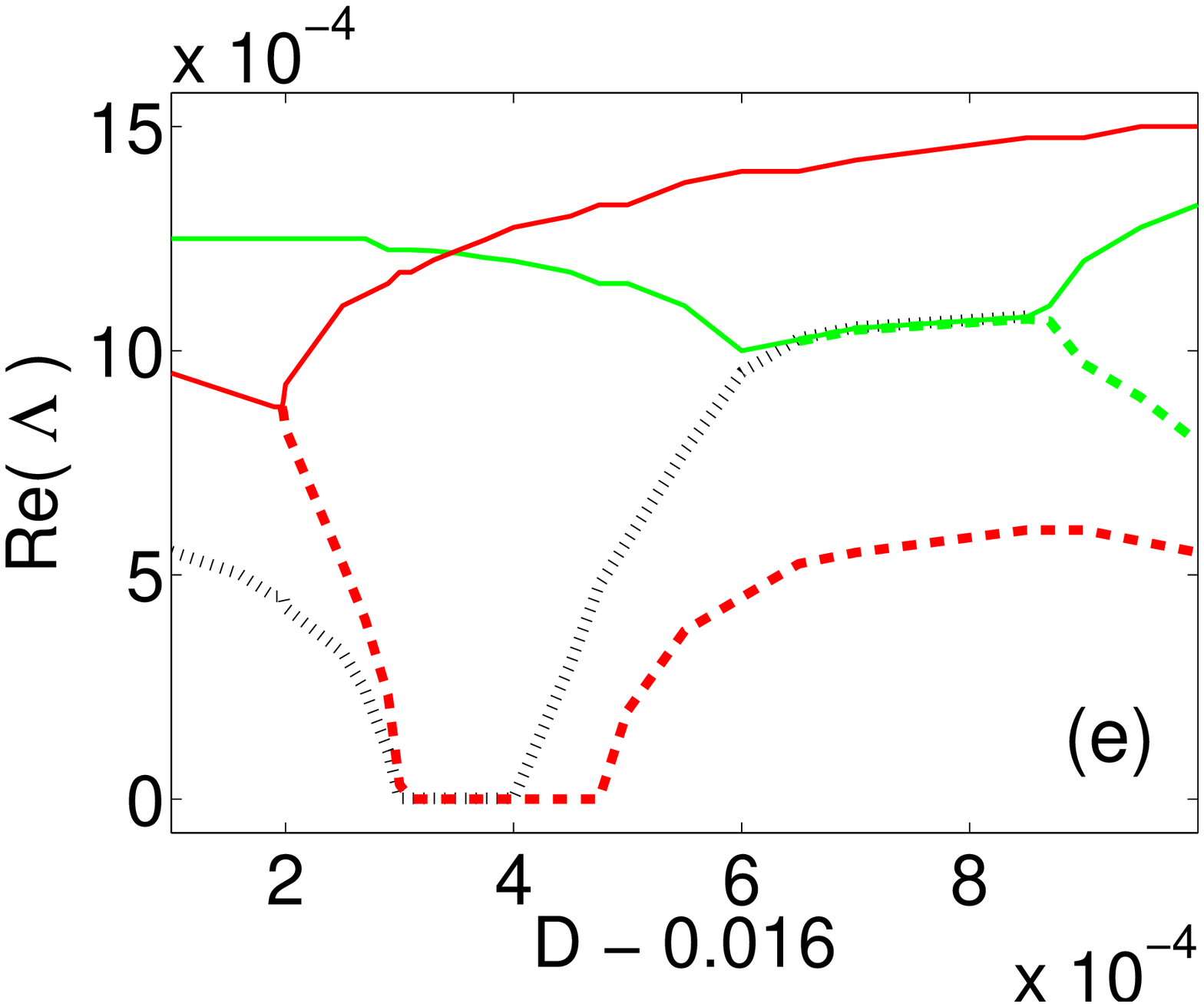}}}
\end{minipage}
 }

\vspace{-2cm}

\hspace*{-0.5cm}
\mbox{ 
\begin{minipage}{5.1cm}
\rotatebox{0}{\resizebox{5.1cm}{7cm}{\includegraphics[0in,0.5in]
 [8in,10.5in]{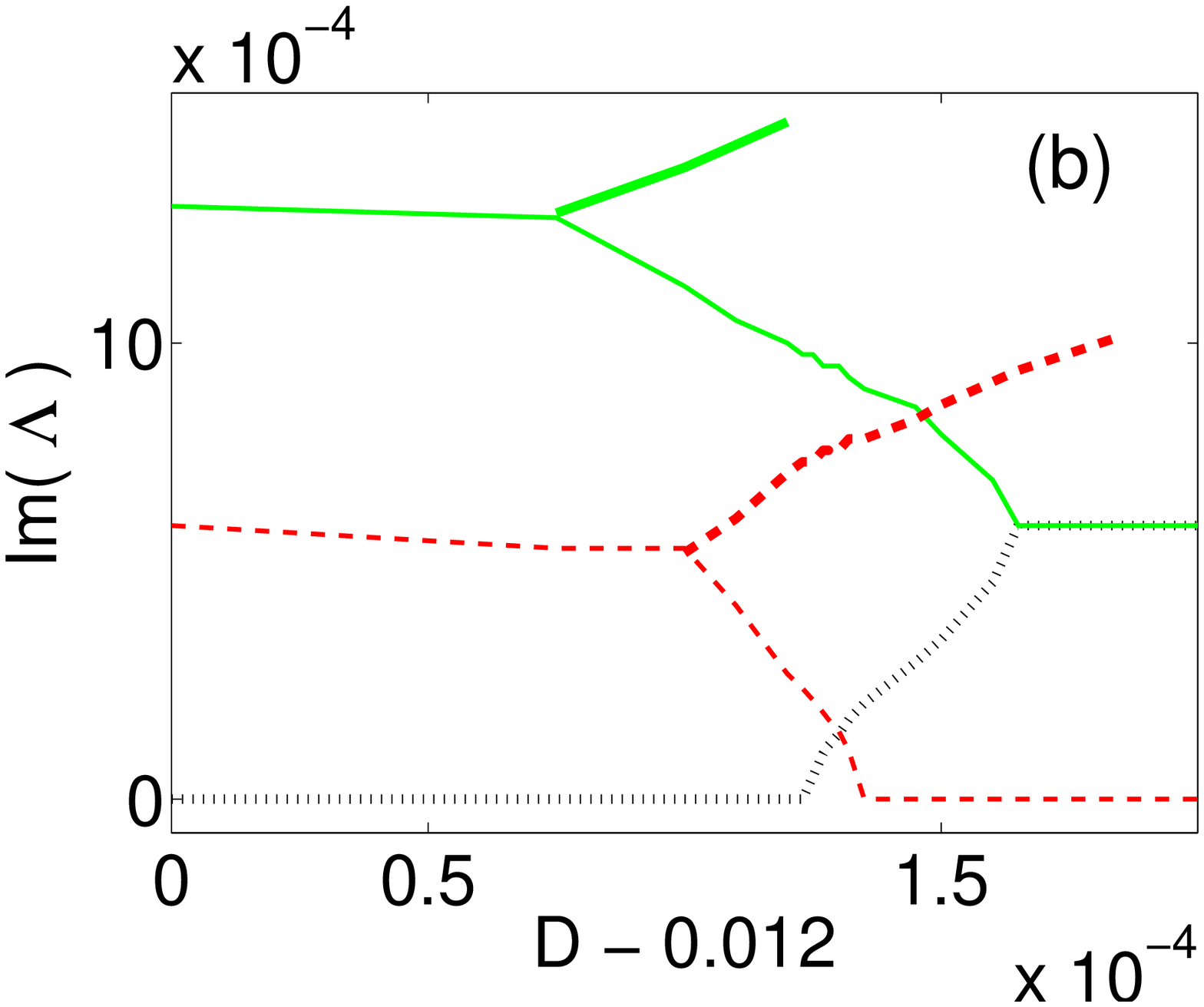}}}
\end{minipage}
\hspace{0.1cm}
\begin{minipage}{5.1cm}
\rotatebox{0}{\resizebox{5.1cm}{7cm}{\includegraphics[0in,0.5in]
 [8in,10.5in]{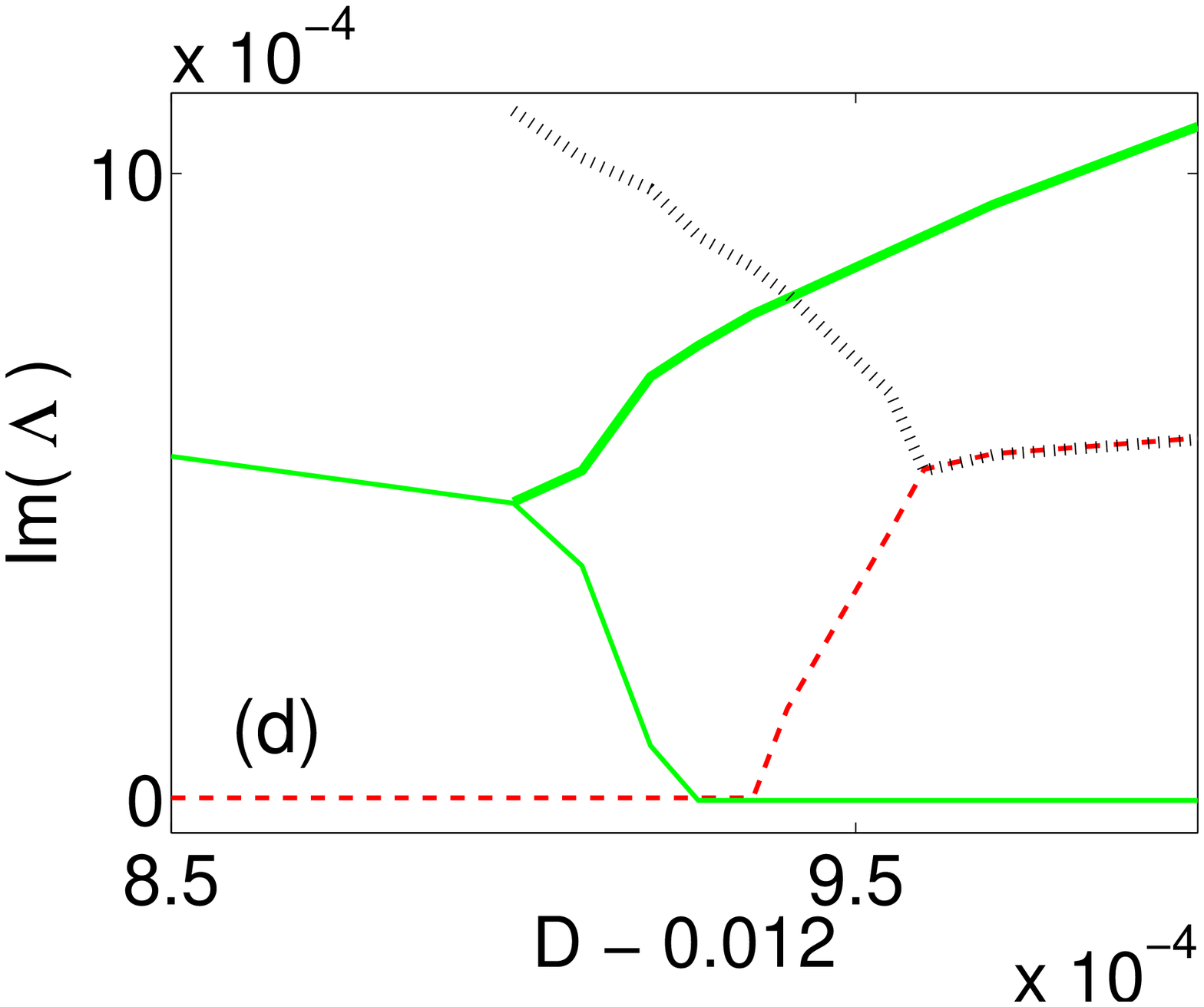}}}
\end{minipage}
\hspace{0.1cm}
\begin{minipage}{5.1cm}
\rotatebox{0}{\resizebox{5.1cm}{7.4cm}{\includegraphics[0in,0.5in]
 [8in,10.5in]{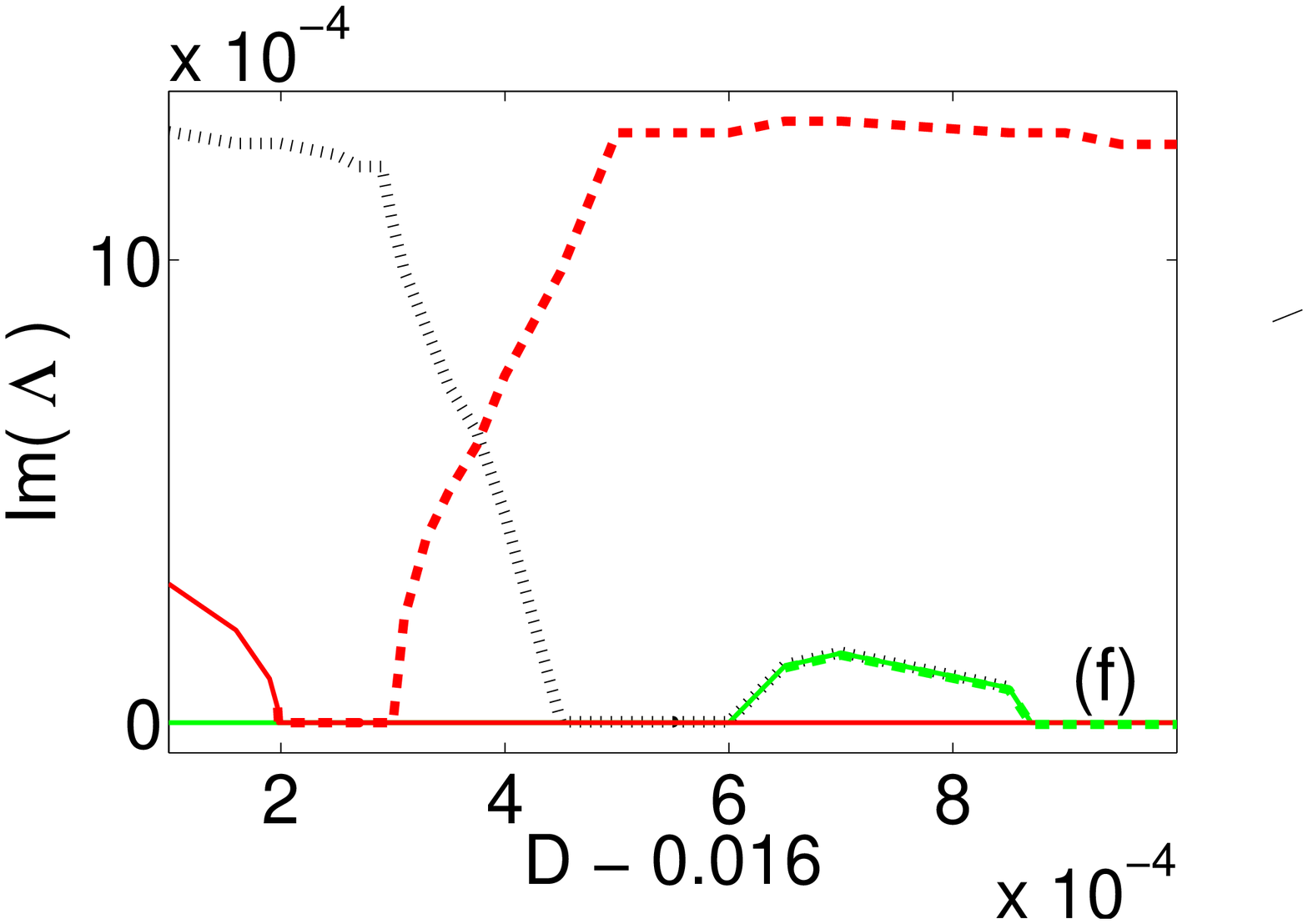}}}
\end{minipage}
 }
\vspace{-1.6cm}
\caption{ (Color online) \ 
Real and imaginary parts of selected modes of \eqref{e_36}, including the
most unstable mode. \ (a) \& (b): $0.01200 \le D \le 0.01220$; \ 
(c) \& (d): $0.01285 \le D \le 0.01300$; \ 
(e) \& (f): $0.01610 \le D \le 0.01700$. Same line colors, styles, and 
widths are used to indicate the same modes within one pair of
panels (e.g., (c) \& (d)). The same line colors/styles/widths in different pairs
of panels (e.g., in (a) \& (b) and (c) \& (d)) do {\em not} imply the same modes.
}
\label{fig_8}
\end{figure}

First, at $D\approx 0.012134$, a ``real" mode appears (see Fig.~\ref{fig_8}(a,b)),
and from this point on there always exists a ``real" mode, even though the particular mode
``born" at $D\approx 0.012134$ disappears later on. Specifically, at $D\approx 0.012928$,
another ``real" mode acquires $\Lambda_R$ greater than that of the mode ``born" at $D\approx 0.012134$,
and the latter mode soon disappears (Fig.~\ref{fig_8}(c,d)). A similar switchover between ``real"
modes occurs at least one more time near $D\approx 0.013750$ (not shown). Next, another ``real"
mode is ``born" via a cascade of bifurcations near $D\approx 0.0162$ (Fig.~\ref{fig_8}(e,f)),
and its $\Lambda_R$ crosses that of the previously dominant-$\Lambda_R$ ``real" mode near
$D\approx 0.01635$. At $D=0.0170$, these two dominant ``real" modes have 
$\Lambda_R\approx 1.4\cdot 10^{-3}$
and $1.5\cdot 10^{-3}$ (Fig.~\ref{fig_8}(e)). Finally, these two modes gradually approach each
other while crossing at least once more near $D=0.01725$. At $D=0.023$ and beyond, their eigenvalues
are the same to five significant figures. Thus, remarkably, the dominant unstable mode 
eventually becomes doubly degenerate. We verified that the same also holds for the
higher-order localized unstable modes.

To conclude, we present a hypothesis as to why the value $C_{\rm cr}$, where
a ``real" mode appears permanently (see above), is near $C=1.012 \,-\, 1.013$. 
Let us interpret \eqref{e_C6} in a way that the $n$ on its r.h.s. is not necessarily
an integer, but
a continuous function of the parameter $D$. For those values of $D$ when $n$ {\em is} an integer,
a mode with a real $\Lambda$ either appears or disappears at the origin $\Lambda=0$. 
Evaluating $n$ at the values
of $D$ listed in the previous paragraph in connection with Figs.~\ref{fig_8}(a)--(d), one finds:
\bsube
\bea
\underline{{\rm at} \;\;D=0.012134}: & & n|_{\nu=1}\approx 29.02, \qquad n|_{\nu=1}-\nu|_{\nu=3}\approx 0.99; 
\label{e_C7a} \\
\underline{{\rm at} \;\;D=0.012928}: & & n|_{\nu=1}\approx30.02, \qquad n|_{\nu=1}-\nu|_{\nu=3}\approx 1.02; 
\label{e_C7b} \\
\underline{{\rm at} \;\;D=0.013750}: & & n|_{\nu=1}\approx31.04, \qquad n|_{\nu=1}-\nu|_{\nu=3}\approx 1.05.
\label{e_C7c} 
\eea
\label{e_C7}
\esube
That is, both $n|_{\nu=1}$ and $n|_{\nu=3}$ are simultaneously very close to integers.
At $D=0.012928$, one of the ``real" modes has not yet disappeared while the next 
one has appeared (Fig.~\ref{fig_8}(c)). From \eqref{e_C7} we observe that at this value of $D$,
the difference $(n|_{\nu=1}-n|_{\nu=3})$ exceeds $1$ for the first time. Thus, we hypothesize
that $C_{\rm cr}\equiv 1+D_{\rm cr}$ 
is found from the condition that $(n|_{\nu=1}-n|_{\nu=3})$ exceeds $1$ 
for the first time. Verification of this hypothesis requires a deeper analytical insight than
we have at the moment. Moreover, finding a value of $C$ past which the dominant real eigenvalue 
increases monotonically (as seen in Fig.~\ref{fig_8}(e)) is also an open question.

\newpage

{\bf\Large List of captions}

\underline{Fig.~1}: \ Growth rate of numerical instability of the s-SSM (a) and
fd-SSM (b) on the plane-wave background. 
The dotted horizontal line indicates how the maximum growth rate depends on the
wave's amplitude. In (a), $k_{m\pi}$, $m=1,2,\ldots$ are the wavenumbers where the
$m$th resonance condition holds (see \cite{ja}): $|\beta|k_{m\pi}^2\dt = m\pi$. 

\medskip

\underline{Fig.~2}: \ See explanation in the text.

\medskip

\underline{Fig.~3}: \ 
Growth rate of the numerical instability for $\dx=40/2^9$ (solid line --- 
analysis of Sec.~IV, stars --- computed by \eqref{e_17}) and for 
$\dx=40/2^{10}$ (dashed line --- 
analysis of Sec.~IV, circles --- computed by \eqref{e_17}).

\medskip

\underline{Fig.~4}: \ 
Normalized phase: \ $|\beta|k^2\dt$ for the s-SSM (dashed) and as given by \eqref{e_19} 
for the fd-SSM (solid). In both cases, $r=5$. The horizontal line indicates the condition
of the first resonance: $|P(k)|=\pi$.

\medskip

\underline{Fig.~5}: \ 
(a): Profiles of the first localized mode on the right side of the soliton
for different values of $C$,
as found by the numerical method of Appendix B. \ (b): Same as in (a), but found 
from the numerical solution of \eqref{e_01}, as explained in the text. \ 
(c): The modes at {\em both} sides of the soliton found from the numerical solution
of \eqref{e_01}. Note that these modes do not ``see" each other because of the
barrier created by the soliton, and hence in general have different amplitudes as
they develop from independent noise seeds. In panels (a)--(c), the potential is 
$\sech^2(\epsilon X)$ (see \eqref{e_35}) and the amplitude of the mode is normalized
to that of the potential. \ (d): Location of the peak of the first localized mode,
found by the method of Appendix B (solid line) and from the solution of \eqref{e_01}
(stars). 
Similar data for $L=40$ and $N=2^{10}$ are very close to those in (d) and hence
are not shown.

\medskip

\underline{Fig.~6}: \ Similar to Fig.~\ref{fig_5}(a), but for the second (a) 
and third (b) localized modes.
 \ (c): $C$ values where localized modes of increasing order appear. Stars --- for 
 $\epsilon = 40/1024$, circles --- for $\epsilon = 40/2048$. 
 
 \medskip
 
\underline{Fig.~7}: \ 
(Color online) \ (a): A solution of \eqref{e_C1} (solid); potential $\sech^2(\epsilon X)$ 
(red dotted). The amplitude of the solution is normalized to that of the potential. \ 
(b): Same as (a), but that panel is ``cut" along the vertical dotted line at the center,
and the resulting halves are interchanged.

\medskip
 
\underline{Fig.~8}: \ 
(Color online) \ 
Real and imaginary parts of selected modes of \eqref{e_36}, including the
most unstable mode. \ (a) \& (b): $0.01200 \le D \le 0.01220$; \ 
(c) \& (d): $0.01285 \le D \le 0.01300$; \ 
(e) \& (f): $0.01610 \le D \le 0.01700$. Same line colors, styles, and 
widths are used to indicate the same modes within one pair of
panels (e.g., (c) \& (d)). The same line colors/styles/widths in different pairs
of panels (e.g., in (a) \& (b) and (c) \& (d)) do {\em not} imply the same modes.

\end{document}